\documentclass{gtart_h}  

\def\ifplaintex{\expandafter\ifx\csname documentclass\endcsname\relax}

\ifplaintex 
\else
\fi

\expandafter\ifx\csname beginpicture\endcsname\relax
\expandafter\ifx\csname documentclass\endcsname\relax
\input pictex \else
\input prepictex \input pictex \input postpictex \fi\fi

\def\gt{{\mathsurround=0pt\it $\cal G\mskip-2mu$eometry \&\ 
$\cal T\!\!$opology}}        

\def\gtp{{\mathsurround=0pt\it $\cal G\mskip-2mu$eometry \&\ 
$\cal T\!\!$opology $\cal P\!$ublications}}  

\def\lognumber#1{\def\thelognumber{#1}}
\def\volumenumber#1{\def\thevolumenumber{#1}}
\def\papernumber#1{\def\thepapernumber{#1}}
\def\volumeyear#1{\def\thevolumeyear{#1}}

\def\pagenumbers#1#2{\def\startpage{#1}\def\finishpage{#2}}
\def\published#1{\def\publishdate{#1}}
\def\proposed#1{\def\theproposer{#1}}
\def\seconded#1{\def\theseconders{#1}}
\def\received#1{\def\receiveddate{#1}}

\def\accepted#1{\def\accepteddate{#1}}

\def\asciiaddress#1{\def\theasciiaddress{#1}}
\def\asciiemail#1{\def\theasciiemail{#1}}

\long\def\asciiabstract#1{\long\def\theasciiabstract{#1}}

\let\\\par\let\thelognumber\relax
\let\thevolumenumber\relax\let\thepapernumber\relax
\let\thevolumeyear\relax\let\thesamplenumber\relax\let\startpage\relax
\let\finishpage\relax\let\publishdate\relax\let\receiveddate\relax
\let\reviseddate\relax\let\accepteddate\relax\let\theasciititle\relax
\let\theasciiauthors\relax\let\theasciiaddress\relax
\let\theasciiabstract\relax
\let\theasciiemail\relax\let\theshortauthors\relax\let\theshorttitle\relax

\long\def\maketitlep{   

\count0=\startpage

\gt\hfill      
\beginpicture
\setcoordinatesystem units <0.33truein, 0.33truein> point at 2.2 0.9
\setplotsymbol ({$\cal G$})
\plotsymbolspacing=9truept
\circulararc 315 degrees from 0 1 center at 0 0
\setplotsymbol ({$\cal T$})
\circulararc 315 degrees from 1 -1 center at 1 0
\endpicture
\break
{\small\ifx\thesamplenumber\relax 
Volume \else Sample
\fi\thevolumenumber\ (\thevolumeyear)
\startpage--\finishpage\nl
Published: \publishdate}
\vglue 0.5truein plus 0.4fil minus 0.1truein

{\parskip=0pt\leftskip 0pt plus 1fil\def\\{\par\smallskip}{\ifplaintex\large
\else\Large\fi\bf\thetitle}\par\medskip}   

\vglue 0pt plus 0.1fil 

{\parskip=0pt\leftskip 0pt plus 1fil\def\\{\par}{\sc\theauthors}
\par\medskip}

\vglue 0pt plus 0.1fil 

{\small\parskip=0pt\let\newline\\
{\leftskip 0pt plus 1fil\def\\{\par}{\sl\theaddress}\par}
\expandafter\ifx\theemail\relax    
\relax\else\vglue 5pt plus 0.02fil minus 2pt\def\\{\stdspace{\rm 
and}\stdspace} 
\cl{Email:\stdspace\tt\theemail}\fi
\ifx\theurl\relax                  
\relax\else\vglue 5pt plus 0.02fil minus 2pt\def\\{\stdspace{\rm 
and}\stdspace}
\cl{URL:\stdspace\tt\theurl}\fi\par}

\vglue 7pt plus 0.3fil minus 3pt

{\bf Abstract}
\vglue 5pt plus 0.1fil minus 2pt

\theabstract

\vglue 7pt plus 0.3fil minus 3pt

{\bf AMS Classification numbers}\quad Primary:\quad \theprimaryclass

Secondary:\quad \thesecondaryclass

\vglue 5pt plus 0.3fil minus 2pt

{\bf Keywords:}\quad \thekeywords

\vglue 10pt plus 0.5fil minus 5pt

{\small  Proposed: \theproposer\hfill Received: \receiveddate\nl
Seconded: \theseconders\hfill 
\ifx\reviseddate\relax                         
Accepted: \accepteddate                        
\else
Revised: \reviseddate                          
\fi}
\eject
}       

\let\maketitlepage\maketitlep
\let\maketitle\maketitlepage

\font\phead=cmsl9 scaled 950
\font=cmsl9 scaled 1050
\font\pnum=cmbx10 scaled 913
\font\lnum=cmbx10 
\font\pfoot=cmsl9 scaled 950
\font=cmsl9 scaled 1050
\ifplaintex
\headline{\vbox to 0pt{\vskip -4.5mm\line{\small\phead\ifnum
\count0=\startpage ISSN 1364-0380 (on line)
1465-3060 (printed) \hfill {\pnum\folio}\else\ifodd\count0\def\\{ }%
\ifx\theshorttitle\relax\thetitle\else\theshorttitle\fi\hfill{\pnum\folio}
\else\def\\{ and }{\pnum\folio}\hfill\ifx\theshortauthors\relax\theauthors
\else\theshortauthors\fi\fi\fi}\vss}}
\footline{\vbox to 0pt{\vglue 0mm\line{\small\pfoot\ifnum\count0=\startpage
\copyright\ \gtp\hfill\else
\gt, Volume \thevolumenumber\ (\thevolumeyear)\hfill\fi}\vss
}}
\else
\makeatletter
\def\@oddhead{{\small\lhead\ifnum\count0=\startpage ISSN 1364-0380 (on line)
1465-3060 (printed) \hfill {\lnum\number\count0}\else\ifodd\count0
\def\\{ }\ifx\theshorttitle\relax \thetitle \else\theshorttitle\fi\hfill
{\lnum\number\count0}\else\def\\{ and }{\lnum\number\count0}
\hfill\ifx\theshortauthors\relax 
\theauthors\else\theshortauthors\fi\fi\fi}}\def\@evenhead{\@oddhead}
\def\@oddfoot{\small\lfoot\ifnum\count0=\startpage\copyright\ \gtp\hfill\else
\gt, Volume \thevolumenumber\ (\thevolumeyear)\hfill\fi}
\def\@evenfoot{\@oddfoot}
\makeatother
\fi

\newwrite\gtoutfile
\long\gdef\makeheadfile{  
{\def\\{, }\def\s{ }
\immediate\openout\gtoutfile head.xxx
\immediate\write\gtoutfile{Proxy-for: \ifx\theasciiauthors\relax
\theauthors\else\theasciiauthors\fi\s<\ifx\theasciiemail\relax\theemail\else\theasciiemail\fi>}
\immediate\write\gtoutfile{\noexpand\\}
\immediate\write\gtoutfile{Authors: \ifx\theasciiauthors\relax
\theauthors\else\theasciiauthors\fi}
{\def\\{ }\immediate\write\gtoutfile{Title: \ifx\theasciititle\relax
\thetitle\else\theasciititle\fi}}
\immediate\write\gtoutfile{Subj-class: GT or SG or MG etc}
\immediate\write\gtoutfile{MSC-class: \theprimaryclass\ifx\thesecondaryclass\relax\else, \thesecondaryclass\fi}
\immediate\write\gtoutfile{Journal-ref: Geom. Topol. \thevolumenumber
(\thevolumeyear) \startpage-\finishpage}
\immediate\write\gtoutfile{Comments: Published by Geometry and Topology at}
\immediate\write\gtoutfile{\s\s http://www.maths.warwick.ac.uk/gt/GTVol\thevolumenumber/paper\thepapernumber.abs.html}
\immediate\write\gtoutfile{\noexpand\\}
\immediate\write\gtoutfile{}
\ifx\theasciiabstract\relax
\immediate\write\gtoutfile{\theabstract}\else
\immediate\write\gtoutfile{\theasciiabstract}\fi
\immediate\write\gtoutfile{}
\immediate\write\gtoutfile{\noexpand\\}
\immediate\write\gtoutfile{}
\immediate\closeout\gtoutfile}}  

\def\maketitlepage{\maketitlep\makeheadfile}
\let\maketitle\maketitlepage

\lognumber{637}
\received{2 August 2005}
\volumenumber{9}\papernumber{47}\volumeyear{2005}
\pagenumbers{2079}{2127}   
\published{27 October 2005}
\accepted{25 October 2005}
\proposed{Ronald Stern}
\seconded{Ronald Fintushel, Simon Donaldson}

\usepackage{amscd,verbatim,amsmath,amssymb}
\newcommand{\p}{\partial}
\renewcommand{\a}{\alpha}
\renewcommand{\b}{\beta}
\newcommand{\Z}{\mathbb Z}
\newcommand{\Q}{\mathbb Q}

\renewcommand{\O}{\mathcal O}
\newcommand{\G}{\mathcal G}
\newcommand{\A}{\mathcal A}
\newcommand{\B}{\mathcal B}
\newcommand{\M}{\mathcal M}
\newcommand{\R}{\mathcal R}

\newcommand{\E}{\mathcal E}
\newcommand{\U}{\mathcal U}
\newcommand{\V}{\mathcal V}
\newcommand{\rhobar}{\bar{\rho}}
\newcommand{\im}{\operatorname{im}}
\newcommand{\ep}{\varepsilon}
\newcommand{\ev}{\operatorname{ev}}
\renewcommand{\phi}{\varphi}

\newcommand{\tr}{\operatorname{tr}}
\newcommand{\su}{\mathfrak{su}}
\newcommand{\so}{\mathfrak{so}}

\newcommand{\Id}{\operatorname{Id}}

\newcommand{\ad}{\operatorname{ad}}

\newcommand{\lf}{\,l\mspace{-2mu}f}
\newcommand{\Spin}{\operatorname{Spin}}
\newcommand{\Ad}{\operatorname{Ad}}

\newcommand{\SO}{\operatorname{SO}}
\renewcommand{\u}{\mathfrak u}
\newcommand{\coker}{\operatorname{coker}}
\newcommand{\Hom}{\operatorname{Hom}}

\newcommand{\hol}{\operatorname{hol}}
\newcommand{\sign}{\operatorname{sign}}
\newcommand{\Stab}{\operatorname{Stab}}
\newcommand{\supp}{\operatorname{supp}}
\newcommand{\Span}{\operatorname{Span}}
\newcommand{\diag}{\operatorname{diag}}
\newcommand{\ZZ}{\Z_2 \oplus \Z_2}

\renewcommand{\P}{\mathcal{PR}}
\newcommand{\lppp}{\lambda'''}
\newcommand{\Tqq}{T^4(q,-q)}
\newcommand{\tqq}{T^3(q,-q)}

\newtheorem{theorem}{Theorem}[section]
\newtheorem{lemma}[theorem]{Lemma}
\newtheorem{proposition}[theorem]{Proposition}
\newtheorem{corollary}[theorem]{Corollary}

\theoremstyle{definition}
\newtheorem{remark}[theorem]{Remark}
\newtheorem{definition}[theorem]{Definition}
\begin{document}

\title{Rohlin's invariant and gauge theory III.\\Homology $4$--tori}
\authors{Daniel Ruberman\\Nikolai Saveliev}
\address{Department of Mathematics, MS 050, Brandeis University\\ 
Waltham, MA 02454, USA\\{\rm\qua and\qua}\\Department of Mathematics,
University of Miami\\PO Box 249085, Coral Gables, FL 33124, USA}
\asciiaddress{Department of Mathematics, MS 050, Brandeis University\\ 
Waltham, MA 02454, USA\\and\\Department of Mathematics,
University of Miami\\PO Box 249085, Coral Gables, FL 33124, USA}

\gtemail{\mailto{ruberman@brandeis.edu}{\rm\qua 
and\qua}\mailto{saveliev@math.miami.edu}}
\asciiemail{ruberman@brandeis.edu, saveliev@math.miami.edu}

\begin{abstract}
This is the third in our series of papers relating gauge theoretic invariants
of certain $4$--manifolds with invariants of $3$--manifolds derived
from Rohlin's
theorem. Such relations are well-known in dimension three, starting with
Casson's integral lift of the Rohlin invariant of a homology sphere. We
consider two invariants of a spin $4$--manifold that has the integral homology
of a $4$--torus. The first is a degree zero Donaldson invariant, counting flat
connections on a certain $SO(3)$--bundle. The second, which depends on the
choice of a 1--dimensional cohomology class, is a combination of Rohlin
invariants of a $3$--manifold carrying the dual homology class. We prove that
these invariants, suitably normalized, agree modulo 2, by showing that they
coincide with the quadruple cup product of 1--dimensional cohomology classes.
\end{abstract}
\asciiabstract{%
This is the third in our series of papers relating gauge theoretic
invariants of certain 4-manifolds with invariants of 3-manifolds
derived from Rohlin's theorem. Such relations are well-known in
dimension three, starting with Casson's integral lift of the Rohlin
invariant of a homology sphere. We consider two invariants of a spin
4-manifold that has the integral homology of a 4-torus. The first is a
degree zero Donaldson invariant, counting flat connections on a
certain SO(3)-bundle. The second, which depends on the choice of a
1-dimensional cohomology class, is a combination of Rohlin invariants
of a 3-manifold carrying the dual homology class. We prove that these
invariants, suitably normalized, agree modulo 2, by showing that they
coincide with the quadruple cup product of 1-dimensional cohomology
classes.}

\maketitle

\section{Introduction}
Let $X$ be a closed smooth spin $4$--manifold which has the integral
homology of $T^4 = S^1\times S^1 \times S^1 \times S^1$. Suppose that
there exists a primitive cohomology class $\alpha \in H^1(X;\Z)$ such
that the infinite cyclic covering $\tilde X_{\alpha}$ corresponding
to $\alpha$ has the integral homology of the $3$--torus. Then $X$ is
called a \emph{$\Z[\Z]$--homology $4$--torus}. The intersection form
on the second cohomology of a $\Z[\Z]$--homology $4$--torus $X$ is
always isomorphic to the sum of three copies of the hyperbolic $2$--form,
but the cup-product on the first cohomology of $X$ may vary.

In this paper, we discuss two invariants of $\Z[\Z]$--homology $4$-tori.
The first one is a Rohlin--type invariant $\bar\rho(X,\alpha)$, which
{\it a priori} depends on the choice of a primitive class $\alpha \in
H^1(X;\Z)$ for which $\tilde X_{\alpha}$ has the integral homology of
the $3$--torus. It is defined in terms of an oriented $3$--manifold $M$
embedded in $X$ that is Poincar\'e dual to the class $\alpha$. Each spin
structure on $X$ induces one on $M$ thus giving rise to a corresponding
Rohlin invariant. The sum of these Rohlin invariants over all the spin
structures on $M$ induced from $X$ is called $\bar\rho(X,\alpha)$, see
Section \ref{S:r-cup}. We show that this is a well defined invariant
with values in $\Z/2\Z$. It extends the Rohlin invariant of
$\Z[\Z]$--homology $S^1\times S^3$, described in~\cite{furuta-ohta}
and~\cite{ruberman:ds}.

The other invariant, $\bar\lambda(X,P)$, comes from gauge theory. It
equals one quarter times the signed count of projectively flat (or
equivalently,
projectively anti-self-dual) connections on an admissible
$U(2)$--bundle $P$ over $X$, see Section \ref{S:l-def}.
This invariant can be viewed as either a degree-zero Donaldson polynomial
of $X$, or a four--dimensional analogue of the Casson invariant. Yet another
way to view $\bar\lambda(X,P)$ is as an extension of the Furuta--Ohta
invariant defined for $\Z[\Z]$--homology $S^1\times S^3$, see 
\cite{furuta-ohta}.

\begin{theorem}\label{T:main}
Let $X$ be a $\Z[\Z]$--homology $4$--torus. For any choice of a primitive
cohomology class $\alpha \in H^1(X;\Z)$ such that $\tilde X_{\alpha}$ has
the integral homology of the $3$--torus, and for any admissible
$U(2)$--bundle $P$,
\[
\bar\rho(X,\alpha) \equiv \bar\lambda(X,P) \equiv \det X\pmod 2
\]
where $\det X$ is defined as $(a_1\cup a_2\cup a_3\cup a_4)\,[X]\pmod 2$,
for any choice of basis $a_1$, $a_2$, $a_3$, $a_4 \in H^1
(X;\Z_2)$.
\end{theorem}

Admissible bundles exist over all $\Z[\Z]$--homology $4$--tori, see Section
\ref{S:admissible}. This fact and Theorem \ref{T:main} together imply
that neither $\bar\rho(X,\alpha)$ nor $\bar\lambda(X,P)\pmod 2$ depend on
choices made in their definition (primitive class $\alpha$ and admissible
bundle $P$, respectively).  The integer invariant $\bar\lambda(X,P) $
does depend on the choice of $P$, as we show by examples in
Section~\ref{S:logtransform}.

Our proof of Theorem~\ref{T:main} takes advantage of a natural action of
$H^1(X;\Z_2) = (\Z_2)^4$ on the moduli space of projectively flat
connections. We identify the latter with the space of projective
representations of $\pi_1 Y$ in $SU(2)$, and use this identification to
show that the orbits with four elements are always non--degenerate and
that the number of such orbits equals $\det X \pmod 2$. In the
non--degenerate situation, this completes the proof because there are no
orbits with one or two elements, and the orbits with eight and sixteen
elements contribute zero to $\bar\lambda(X,P)\pmod 2$. The general
case is handled similarly after one finds a generic perturbation which
is $H^1(X;\Z_2)$--equivariant.

Because of the equivariance condition, constructing such a
perturbation requires a rather delicate argument, which is carried
out in Sections 6, 7, and 8.   The different orbits (of size $4$,
$8$, or $16$) are handled separately, with no perturbation required
at the $4$-orbits.   The hardest case is that of the $8$-orbits, in
which connections have a $\Z_2$ stabilizer inside $H^1(X;\Z_2)$.   We
consider first a perturbation that affects the deformation complex
tangential to an orbit, and for this we perturb the ASD equations
using holonomy functions, as
in~\cite{floer:instanton,taubes:casson,donaldson:orientation,herald:perturbations,kronheimer:higher-rank}.

After this perturbation, the orbit is $0$-dimensional, but the points
in it may not be smooth points in the moduli space as there is
potentially some nontrivial cohomology in the deformation complex
normal to the orbit.   If so, we make a further perturbation using
the abstract perturbations found in~\cite{donaldson}
and~\cite{furuta:perturbation}.   Having achieved non-degeneracy at 
the 8--orbits
using equivariant
perturbations, we move on to the non-degeneracy at the 16--orbits in
Section~\ref{S:sixteen}.

In \cite{ruberman-saveliev:casson}, we studied a three--dimensional analogue
of the invariant $\bar\lambda(X,P)$, denoted $\lambda'''(Y,w)$. It is
obtained by counting projectively flat connections on a $U(2)$--bundle $P$
over a homology $3$--torus $Y$ with $c_1(P) = w \ne 0\pmod 2$. We showed
that $\lambda'''(Y,w) \equiv \det Y \pmod 2$. Conjecturally,
the invariant $\lambda'''(Y,w)$ coincides with the Lescop invariant
\cite{lescop:casson} obtained via a surgery approach. Lescop's invariant
is independent of the choice of $w$ and in fact equals $(\det Y)^2$. It is
natural to wonder if similar results hold for $\bar\lambda(X,P)$; at the
end of the paper we present a family of examples that show that
$\bar\lambda(X,P)$ is not determined by $\det X$, and that in fact it
depends on the choice of bundle $P$.

Theorem \ref{T:main} is part of a broader program for relating
gauge--theoretic invariants of some simple $4$--manifolds to Rohlin-type
invariants~\cite{ruberman-saveliev:casson,ruberman-saveliev:mappingtori}.
A survey of this work appears in~\cite{ruberman-saveliev:survey}.
In the spirit of Witten's conjecture \cite{witten:monopole} relating 
Donaldson and
Seiberg--Witten invariants, Theorem \ref{T:main} is also consistent with
the result of the first author and S.~Strle \cite{ruberman-strle:tori} on
the (mod 2) evaluation of the Seiberg-Witten invariant for homology
$4$-tori.

We would like to express our appreciation to Chris Herald for
pointing out a gap in our treatment of equivariant transversality in
an earlier version of this paper, and for an extensive correspondence
on the subject.  We would also like to thank Cliff Taubes for
suggesting the strategy used in Section~\ref{S:tau} for achieving
transversality at the orbits of $H^1(X;\Z_2)$ with eight elements.
The first author was partially supported by NSF Grants
9971802 and 0204386. The second author was partially supported by NSF
Grant 0305946.


\section{Algebraic topology of homology tori}\label{S:algtop}
The discussion of both the Rohlin and Donaldson invariants requires some
information about algebraic topology of homology tori;  relevant results
are collected in this section.

An \emph{$n$--dimensional homology torus} $X$ is a manifold of dimension
$n$ such that $H_* (X;\Z) = H_* (T^n;\Z)$ where $T^n = S^1\times\ldots
\times S^1$ ($n$ times). Let $a_1,\ldots, a_n$ a basis in $H^1 (X;\Z_2)$
then $\det X = (a_1\cup\ldots\cup a_n)\,[X]\pmod 2$ is independent of
the choice of $a_1,\ldots, a_n$ and is called the \emph{determinant} of
$X$. A homology torus $X$ is called \emph{odd} if $\det X = 1\pmod 2$;
otherwise, $X$ is called \emph{even}.

Let $R$ be a commutative ring with unity; at various points we will use
the rings $\Z$, $\Q$, $\Z_2$, and the ring $\Z_{(2)}$ of integers localized
at $2$, i.e. with all odd primes inverted. An \emph{$n$--dimensional
$R$--cohomology torus} is a manifold $X$ such that $H^*(X;R)$ and
$H^*(T^n;R)$ are isomorphic as rings.

\begin{lemma}\label{L:cup}
Suppose that $X$ is an $n$--dimensional odd homology torus. Then $X$ is
an $n$--dimensional $\Z_2$--cohomology torus and a $\Q$--cohomology torus.
\end{lemma}

\begin{proof} By the universal coefficient theorem,
the dimension of $H^k(X;\Z_2)$ is the same as the rank of
$H^k(X;\Z)$, and a basis $a_1,\ldots, a_n$ for $H^1(X;\Z)$ gives
rise to a basis for $H^1(X;\Z_2)$.  It is straightforward to show
that the $\binom{n}{k}$ distinct cup products of the $a$'s are
linearly independent in $H^k(X;\Z_2)$, and thus give a basis.  The
proof with $\Q$ coefficients is the same.
\end{proof}

We remind the reader that, for any $Y$, the cup-square $H^2(Y;\Z_2) \to
H^4(Y;\Z_2)$ actually lifts to a map $H^2(Y;\Z_2) \to H^4(Y;\Z_4)$ known
as the Pontrjagin square~\cite{pontrjagin:pi3}.  In our situation, where
every $\Z_2$ cohomology class lifts to an integral class, the Pontrjagin
square can be computed by lifting to an integral class, and reducing the
cup product mod $4$.  We will write this with the usual notation for cup
product.

\begin{corollary}\label{C:cup}
Let $X$ be an $n$--dimensional odd homology torus then the cup-product
\begin{equation}\label{E:cup}
\cup: \Lambda^2 H^1 (X;\Z_2) \to H^2 (X;\Z_2)
\end{equation}
gives rise to a bijective correspondence between decomposable elements
in $\Lambda^2 H^1(X;\Z_2)$ and elements $w\in H^2 (X;\Z_2)$ such that
$w\cup w = 0\pmod 4$.
\end{corollary}

\begin{proof}
This is true for $T^n$ hence also for $X$, since the $\Z_2$--cohomology
rings of $X$ and $T^n$ are isomorphic.
\end{proof}

The natural map $X \to K(\pi_1 X, 1)$ induces a monomorphism (see
\cite{brown:cohomology})
\begin{equation}\label{E:iota}
\iota: H^2 (\pi_1 X;\Z_2)\to H^2 (X;\Z_2).
\end{equation}

\begin{corollary}\label{C:odd}
Let $X$ be an odd homology torus then the map $\iota$ is an isomorphism.
\end{corollary}

\begin{proof}
This follows from the commutative diagram
\[
\begin{CD}
\Lambda^2 H^1 (\pi_1 X; \Z_2) @>\cong >> \Lambda^2 H^1 (X;\Z_2) \\
@V\cup VV    @V\cup VV \\
H^2 (\pi_1 X; \Z_2) @>\iota >> H^2 (X;\Z_2)
\end{CD}
\]
whose upper arrow is an isomorphism because $H^1 (\pi_1 X; \Z_2) = H^1 (X;
\Z_2)$, and the right arrow is an isomorphism by Lemma \ref{L:cup}. Since
$\iota$ is injective, the remaining two arrows in the diagram are also
isomorphisms.
\end{proof}

\begin{lemma}\label{L:rings}
A manifold $X$ is a $\Z_{(2)}$-cohomology torus if and only if it is a
$\Z_2$-cohomology and $\Q$-cohomology torus.
\end{lemma}

\begin{proof}
The statement about the cohomology groups is  a simple application of
the universal coefficient theorem and Poincar\'e duality.  To
understand the ring structure, note that $H^1(X;\Z_{(2)}) \cong
H^1(X;\Z) \otimes \Z_{(2)}$.  Hence the $n$-fold cup product on
$H^1(X;\Z_{(2)})$ is odd if and only if the $n$-fold cup products on
$H^1(X;\Q)$ and $H^1(X;\Z_2)$ are nontrivial.
\end{proof}

In this paper, we will be concerned with  the following situation: we are given
a homology $4$--torus $X$, together with a specific double covering.
This double covering corresponds to a cohomology class $\alpha \in
H^1(X;\Z_2)$, and will be denoted $X_{\alpha} \to X$.

\begin{proposition}\label{P:toruscover} Suppose that $X$ is a
$\Z_{(2)}$-cohomology torus, and that $\pi: X_{\alpha}\to X $ is a
nontrivial double covering corresponding to $\alpha\in H^1(X;\Z_2)$.
Then $X_\alpha$ is a $\Z_{(2)}$-cohomology torus.
\end{proposition}

\begin{proof}
Associated to the covering $\pi:X_\alpha \to X$ is the Gysin sequence in
cohomology (where the coefficients are understood to be $\Z_2$)
\[
\begin{CD}
\to H^{k-1}(X)  @>\cup \alpha>> H^k(X) @>\pi^*>> H^k(X_\alpha) @>>>
H^k(X) @>\cup \alpha>> H^{k+1}(X)\\
@V\cong VV @V\cong VV @V= VV @V\cong VV @V\cong VV \\
\to \Z_2^{\binom{n}{k-1}}  @>\cup \alpha>> \Z_2^{\binom{n}{k}}  @>\pi^*>>
H^k(X_\alpha) @>>> \Z_2^{\binom{n}{k}} @>\cup \alpha>> \Z_2^{\binom{n}{k+1}}
\end{CD}
\]
By Lemma~\ref{L:rings}, the $\Z_2$-cohomology ring of $X$ is an exterior
algebra over $\Z_2$. If we choose a generating set $\alpha = \alpha_1$,
$\alpha_2,\ldots,\alpha_n$, then the image of the map $\cup\,\alpha$ in
$H^k(X;\Z_2)$ is the $\binom{n-1}{k-1}$ dimensional subspace spanned by
products of $\alpha$ with monomials not involving $\alpha$. Likewise, the
kernel of $\cup\,\alpha: H^k(X;\Z_2) \rightarrow H^{k+1}(X;\Z_2)$ has
dimension $\binom{n-1}{k-1}$. From exactness, the dimension of
$H^k(X_\alpha;\Z_2)$ is $\binom{n}{k}$. It is straightforward to check, by
similar arguments, that the $n$-fold cup product on $H^1(X_\alpha;\Z_2)$
is non-trivial.
\end{proof}

\medskip

Finally, we record for later use some facts about the relation between the
cohomology of a $\Z[\Z]$--homology $4$--torus $X$ and that of its infinite
cyclic cover $\pi: \tilde X_\alpha \to X$ classified by $\alpha \in H^1
(X; \Z)$.

\begin{lemma}\label{L:zcover}
Let $X$ be a $\Z[\Z]$--homology $4$--torus whose infinite cyclic cover
$\tilde X_\alpha$ has the integral homology of the $3$--torus.  Then
\begin{enumerate}
\item
The cup-product pairing $H^1 (\tilde X_\alpha;\Z_2) \otimes
H^2(\tilde X_\alpha;\Z_2) \to H^3(\tilde X_\alpha; \Z_2) \cong \Z_2$ is
non--degenerate.
\item $\pi^*: H^i (X;\Z_2) \to H^i (\tilde X_\alpha;\Z_2)$ is a surjection
for all $i$.
\item
$\ker \pi^*:H^1 (X;\Z_2) \to H^1 (\tilde X_\alpha;\Z_2)$ is the subgroup
generated by $\alpha$.
\end{enumerate}
\end{lemma}

\begin{proof}
The universal coefficient theorem implies that the  $\Z_2$ homology of
$\tilde X_\alpha $ is finitely generated, because the integral homology
is finitely generated.  Hence the first item  follows from Milnor's
duality theorem~\cite{milnor:covering}. The other two parts are proved
using the long exact sequence (where the coefficients
are understood to be $\Z_2$)
\begin{equation}\label{E:coverseq}
\begin{CD}
@>>>  H^i(X) @>\pi^*>> H^i(\tilde X_\alpha)  @> t^* -1>>
H^i(\tilde X_\alpha) @>\delta>> H^{i+1}(X)  @>\pi^*>>
\end{CD}
\end{equation}
derived in~\cite[proof of assertion 5]{milnor:covering}. Here $t$
represents the covering translation.  Looking at the beginning of
this sequence, we see that $t^*: H^0 (\tilde X_\alpha;\Z_2) \to
H^0(\tilde X_\alpha;\Z_2)$ is the identity, since $\tilde X_\alpha$
is connected. Thus the image of $\delta: H^0(\tilde X_\alpha;\Z_2)
\to H^1(X;\Z_2)$ is one--dimensional.  The rest follows by counting
dimensions and using exactness.
\end{proof}
We remark that the same proof would work with any field coefficients,
from which we can deduce that the conclusions of Lemma~\ref{L:zcover}
hold with integral coefficients as well.  In the proof of
Corollary~\ref{C:cup12} below, we need the following claim.

\begin{lemma}
In the exact sequence~\eqref{E:coverseq} with integral coefficients,
we have $(\delta c)^2 =0$ for any $c\in H^1(\tilde X_\alpha;\Z)$.
\end{lemma}

\begin{proof}
Consider the following commutative square:
$$
\begin{CD}
H^1(\tilde X_\alpha;\Z) @>\delta>> H^{2}(X;\Z)   \\
@VV\cap [M] V   @VV\cap [X] V\\
H_2(\tilde X_\alpha;\Z) @>\pi_*>> H_{2}(X;\Z)
\end{CD}
$$
where the right-hand vertical arrow is Poincar\'e duality, and
$M^3$ is a lift to
$\tilde X_\alpha$ of the Poincar\'e dual of $\alpha$.  Now
$\langle
(\delta c)^2, [X]\rangle = (\delta c \cap [X])^2$
by the
correspondence between intersections and cup
products.  Since $(\delta c \cap [X])^2 = (\pi_*(c \cap [M]))^2$, it
suffices to show the vanishing of intersection numbers between
classes in $H_2(X;\Z)$ lying in the image of $\pi_*$. But the map $M
\to \tilde X_\alpha$ has degree one, so the induced map $H_2 (M;\Z)
\to H_2(\tilde X_\alpha;\Z)$ is onto.  So if $a$, $b \in H_2(\tilde
X_\alpha;\Z)$, we can represent both classes by cycles in $M$, and so
$\pi_*(a)$ and $\pi_*(b)$ are represented by cycles lying in $M
\subset X$.  It follows readily that $\pi_*(a)\cdot \pi_*(b)= 0$.
\end{proof}

Using this, we get the following corollary, which will be used in
Section~\ref{S:admissible} to construct admissible bundles over
$\Z[\Z]$--homology $4$-tori.

\begin{corollary}\label{C:cup12}
Let $X$ be a $\Z[\Z]$--homology $4$--torus.  Then there are classes $w\in
H^2(X;\Z_2)$ and $\xi \in H^1(X;\Z_2)$ with $w \cup w \equiv 0 \pmod{4}$
and $w\cup \xi \neq 0$.
\end{corollary}

\begin{proof}
We will show that there is a class $w\in H^2(X;\Z_2)$ with $w\cup w\equiv
0 \pmod{4}$ such that $\pi^*(w) \neq 0 \in H^2(\tilde X_\alpha;\Z_2)$. By
Milnor duality on $\tilde X_\alpha$, see Lemma \ref{L:zcover}, there is a
class $\xi' \in H^1 (\tilde X_\alpha;\Z_2)$ with $\pi^*(w) \cup \xi' \neq
0$. But $\xi' = \pi^*\xi$ for some $\xi \in H^1(X;\Z_2)$, see Lemma
\ref{L:zcover}, and therefore $w\cup \xi \neq 0$ by naturality of the cup
product.

To show the existence of such a class $w$, consider the exact
sequence~\eqref{E:coverseq}, which shows that $H^2(X;\Z_2)\cong \im(\delta)
\oplus V$ where both summands are $3$--dimensional vector spaces over $\Z_2$
and $\pi^*:V \rightarrow H^2(\tilde X_\alpha; \Z_2)$ is an isomorphism.
Suppose that there is no element $0\neq a \in V$ with $a \cup a \equiv 0
\pmod 4$. Since $X$ is spin, we then have that $a \cup a = 2 \pmod 4$, so
the formula
\begin{equation}\label{E:square}
(a+b)^2 \equiv a^2 + 2\,a\cup b +b^2 \pmod{4}
\end{equation}
implies that $a\cup b \equiv 1\pmod 2$ for all $a\neq b\in V$. It is easy to
see (say, by considering a basis for $V$) that this makes the cup-product
form on $V$ singular. By non-singularity of the cup product on $X$, there
is a non-zero cup product $a\cup b$ where $a\in V$ and $b\in \im(\delta)$.
By the previous lemma, $b^{\,2} \equiv 0\pmod 4$, since $b$ lifts to an
integral class whose square is $0$. Applying equation~\eqref{E:square} once
more, we conclude that $(a+b)^2 \equiv 0 \pmod 4$; by construction $\pi^*
(a+b) = \pi^*(a) \neq 0$. To finish the proof, we set $w = a + b$.
\end{proof}


\section{The Rohlin invariants}
After describing a Poincar\'e duality theorem for smooth periodic-end
manifolds, we introduce a Rohlin-type invariant for homology $4$-tori,
and then show that a certain combination of such Rohlin invariants is
determined by the cup-product structure.

\subsection{Periodic-end manifolds}\label{Sec:periodic}
A periodic-end manifold $W_\infty$, in the sense of Tau\-bes
\cite{taubes:periodic}, has the following structure. The model for the end
is an infinite cyclic covering $\pi:
\tilde X\to X$ classified by a cohomology class $\alpha\in H^1(X;\Z)$. This
choice implies a specific generator $t$ of the covering translations and
thereby picks out one end of $\tilde X$ as lying in the positive direction.
We denote this end by $\tilde X_+$. Choose a codimension-one submanifold
$M \subset X$ dual to $\alpha$, and a lift $M_0$ of $M$ to  $\tilde X$,
giving copies $M_k = t^k(M_0)$. Thus $\tilde X$ is decomposed into copies
$X_k = t^k(X_0)$ of a fundamental domain $X_0$, which is a manifold with
boundary $-M_0 \cup M_1$. Write
$$
N_k = \bigcup_{m\ge k} X_m.
$$
Finally, suppose we have a compact manifold $W$ with boundary $M$.
The end-periodic manifold $W_\infty$ is given by $W_\infty = W
\cup_{M_0} N_0$.  In other words, $W_\infty$ is a non-compact
manifold, with a single end that coincides with the positive end of
an infinite cyclic cover.

The duality theorem is essentially a combination of duality theorems
of Milnor~\cite{milnor:covering} and Laitinen~\cite{laitinen:ends};
we briefly review these results.  For a non-compact manifold $Y$,
define the \emph{end cohomology} $H^j_e(Y)$
to be the direct limit of the groups $H^i(Y-K;\Z)$ where $K$ runs over
an exhausting set of compact sets.  It is not hard to show that
$H^j_e(Y)$ splits as a direct product of groups $H^j_e(Y_\epsilon)$
where $Y_\epsilon$ runs over the set of ends of $Y$.  Dually,
Laitinen defines the end homology in terms of (the inverse limit of)
the algebraic mapping cone of $C_*(Y-K)\to C_*(Y)$.  Again, it is
not hard to define the end homology of a single end.  It is evident
from the definition that the end cohomology depends only on the end,
or in other words that two manifolds that coincide off a compact set
have the same end cohomology.  The same is true about the end
homology, although this requires more work.    Laitinen proved that
cap product with an appropriate fundamental class gives
Poincar\'e-Lefschetz duality; in particular the end cohomology in
dimension $k$ is isomorphic to the end homology in dimension $n- k
-1$, where $n = \dim X$.

Recall that Milnor~\cite{milnor:covering} showed that if $X$ is
oriented and $\tilde X$ has finitely generated homology with
coefficients $F$ (for $F$ a field or the integers), then $\tilde X$
satisfies Poincar\'e duality:
\[
\cap\, [M]: H^j(\tilde X; F)\stackrel{\cong}{\rightarrow}
H_{n-j-1}(\tilde X;F).
\]
In other words, $\tilde X$ resembles a compact manifold of dimension
$n-1$. Our duality theorem is a relative version of this.

\begin{theorem}\label{T:poincare}
Let $W_\infty$ be an oriented periodic end-manifold, as described
above.  Assume that the homology of the infinite cyclic cover $\tilde
X$, with coefficients in $F = $ a field  or the integers, is finitely
generated.   Then the groups $H^i(W_\infty;F)$ and $H_i(W_\infty;F)$
are  finitely generated, and there are long exact sequences related
by Poincar\'e duality isomorphisms\eject
\begin{equation}\label{E:les}
\minCDarrowwidth1pc
\begin{CD}
\to H^j_c(W_\infty;F)  @>>> H^j (W_\infty;F)  @>>> H^j(\tilde X;F)
@>>> H^{j+1}_c(W_\infty;F)\to \\
@VV\cap[W_\infty]V @VV\cap[W_\infty]V @VV\cap[M]V @VV\cap[W_\infty]V \\
\to H_{n-j}(W_\infty;F)  @>>> H^{\lf}_{n-j} (W_\infty;F)  @>>>
H_{n-j-1}(\tilde X;F)  @>>> H_{n-j-1}(W_\infty;F)\to
\end{CD}
\minCDarrowwidth2.5pc
\end{equation}
where $H^*_c$ stands for cohomology with compact support, and $H_*^{\lf}$
for locally finite homology.
\end{theorem}

\noindent
In other words, $W_\infty$ behaves like an $n$--manifold with boundary,
where the homological `boundary' $\tilde X$ looks like an
$(n-1)$--manifold.

\begin{proof}
Laitinen noted that there is an
exact sequence like~\eqref{E:les}, with $H^j_e (Y)$ in place of
$H^j(\tilde X)$, and proved that the end homology and cohomology
groups satisfy Poincar\'e duality of dimension $n-1$.   The end
homology sits in an exact sequence involving the ordinary homology
and the locally finite homology of $Y$.  There are related
definitions in Milnor~\cite[page 124 and footnote]{milnor:covering} of
the homology/cohomology of a space relative to a single end.

To prove the theorem, we show that the finitely-generated hypothesis
implies that the end homology and cohomology of $W_\infty$ coincide
with the ordinary homology and cohomology of $\tilde X$.  Note that
the end cohomology and homology of $W_\infty$ are the same as the end
cohomology and homology of the positive end $\tilde X_+$ of $\tilde X$.
Let $N_k$ be a neighborhood of the positive end of $\tilde X$, so there
is a long exact sequence
\begin{equation}\label{E:xles}
                  \cdots \to H^j(\tilde X, N_k) \to H^j(\tilde X) \to H^j(N_k)
\to \cdots
\end{equation}
Passing to the limit as $k\to \infty$, we get a long exact sequence
\begin{equation}\label{E:endles}
                  \cdots \to H^j(\tilde X, \tilde X_+) \to H^j(\tilde X) \to
H^j(\tilde X_+) \to \cdots
\end{equation}
But~\cite[Assertion 8]{milnor:covering} says that under our
hypotheses, the relative term vanishes, and so we get the desired
isomorphism.

The statement about homology ends is proved similarly, using the
analogue of the exact sequence~\cite[(2.6)]{laitinen:ends} for the
homology of a single end of $\tilde X$.
\end{proof}

For the rest of this subsection, assume that $W_\infty$ is an
oriented end-periodic $4$--manifold, such that the end cohomology with
rational coefficients is finitely generated.  Note that the cup
product $H^2_c(W_\infty;\Q) \otimes H^2_c(W_\infty;\Q) \to
H^4_c(W_\infty;\Q) \cong \Q$ is well-defined, and hence has a
signature which we denote by $\sign(W_\infty)$.  Now $W_\infty = W
\cup_{M_0} N_0$, and it is easy to see that $N_0$ also has a
signature.   Any of these signatures can be computed using
intersections rather than cup products.

\begin{lemma}\label{L:signature}
If the end cohomology is finitely generated, then $\sign(W_\infty) =
\sign(W)$.
\end{lemma}

\begin{proof}
The usual proof of Novikov additivity applies to show that
$\sign(W_\infty) = \sign(W) +\sign(N_0)$, so it suffices to show that
$\sign(N_0) = 0$.  Now $M_0$ carries the fundamental class of $\tilde
X$, so that the map $H_2(M_0;\Q) \to H_2(\tilde X;\Q)$ is surjective.
Therefore, there are no non-trivial intersections between cycles in
$\tilde X$, and so the same is true for $N_0$.  It follows that
$\sign(N_0) = 0$.
\end{proof}


\subsection{The Rohlin invariants}
For any manifold $Y$, denote by $\Spin(Y)$ the set of $\Spin$ structures
on $Y$. Let $X$ be a homology $4$--torus, that is, a smooth spin $4$--manifold
with $H_*(X;\Z) = H_*(T^4;\Z)$. Choose a spin structure $\sigma \in
\Spin(X)$ and a primitive cohomology class $\alpha\in H^1(X;\Z)$. These
ingredients are enough to define an invariant $\rho(X,\alpha,\sigma) \in
\Q/2\Z$ if $X$ is odd.  On the other hand, proving well-definedness of
$\rho(X,\alpha,\sigma)$ for even $X$ requires an additional assumption
about the infinite cyclic covering $\tilde X_\alpha\to X$ corresponding
to $\alpha$. Here is the definition; the fact that it is independent of
choices is the main result of this subsection.

\begin{definition}\label{D:rhodef}
Let $(X,\alpha,\sigma)$ be as above, and let $M$ be a smooth oriented
$3$--manifold embedded in $X$ that is Poincar\'e dual to $\alpha$. Then
$\sigma$ induces in a canonical way a spin structure on $M$ which we
again call $\sigma$, and we define
$$
\rho(X,\alpha,\sigma) = \rho(M,\sigma) \in \Q/2\Z
$$
where $\rho(M,\sigma)$ is the usual Rohlin invariant of $M$.  Our
convention is that $\rho(M,\sigma) = 1/8\,\sign(V)$, where $V$ is any
smooth compact spin manifold with boundary $M$.
\end{definition}

\begin{theorem}\label{T:rhowelldef}
If $X$ is an odd homology $4$--torus, then $\rho(X,\alpha,\sigma)$ is
independent of the choice of submanifold $M$.  The same is true if
$X$ is an even homology $4$--torus, and the rational homology of
$\tilde X_\alpha$ is finitely generated.  In either case, if the
integral homology  of $\tilde X_\alpha$ is the same as the homology
of the $3$--torus, then $\rho(X,\alpha,\sigma) \in \Z/2\Z$.
\end{theorem}

\begin{proof}
Suppose first that $X$ is a homology torus for which the rational homology of
$\tilde X_\alpha$ is finitely generated.  It follows from Milnor's duality
theorem~\cite{milnor:covering} that $\tilde X_\alpha$ satisfies Poincar\'e
duality over $\Q$ as if it were a $3$--manifold.  Note that any $M \subset X$
dual to $\alpha$ lifts to $X_\alpha$, and gives rise to a fundamental class
for $\tilde X_\alpha$. Thus the inclusion $M \to \tilde X_\alpha$ is a map of
non-zero degree, and thus is surjective in rational homology. In particular,
any classes $x,y \in H_2(\tilde X_\alpha;\Q)$ have trivial intersection number.

Now suppose that $M, M'$ are two submanifolds dual to $\alpha$.  It is easy to
show that they have disjoint lifts to $\tilde X_\alpha$, and so cobound a spin
submanifold $V \subset \tilde X_\alpha$.  The intersection form on $V$ is
trivial, by the previous paragraph, so its signature is trivial.  Hence the
Rohlin invariants of $M$ and $M'$ must agree.

Now suppose that $\tilde X_\alpha$ has the integral homology of the
$3$--torus, and let $M^3$ be an oriented submanifold of $X$ dual to
$\alpha$ with a spin structure $\sigma$ induced from $X$. Choose a spin
$4$--manifold $W$ with boundary $M$ and a lift $M_0$ of $M$ to $\tilde
X_\alpha$.   Form, as in Section~\ref{Sec:periodic}, a periodic
end-manifold $W_\infty$ whose end is modeled on the positive end of $\tilde
X_\alpha$. By definition, $\rho(M,\sigma)$ equals one eighth of
$\sign(W)$. The latter signature coincides with
the signature of $W_\infty$ by Lemma~\ref{L:signature}.  Recall that
for a compact $4$--manifold whose boundary has no torsion in $H_1$,
the intersection form splits as a unimodular form plus the radical.
This follows from Poincar\'e duality, including the relation between
duality on the manifold and on its boundary (as expressed in the
exact sequence~\eqref{E:les}).    As we showed in
Theorem~\ref{T:poincare}, the same duality holds for a periodic end
manifold.   Hence with the hypothesis on $\tilde X_\alpha$, the
intersection form on $W_\infty$ is a sum of a unimodular form and the
trivial form.  But (cf.~\cite{milnor-husemoller}) the signature of an
even unimodular form is always divisible by $8$, and so the Rohlin
invariant lies in $\Z/2\Z$.

Finally, we verify that the Rohlin invariant of an odd homology torus
is well-defined, without any additional hypothesis on the infinite
cyclic cover. The main observation is the following. Suppose that $X$
is a $\Z_{(2)}$ cohomology torus, and that $M$ is dual to a
non-trivial class $\alpha \in H^1(X;\Z)$. Then the intersection form
on $X - M$ is trivial.  To see this, suppose that $\Sigma$ is a
surface in $X$, disjoint from $M$, with $\Sigma \cdot \Sigma \neq 0$.
Setting $\beta = \mathrm{PD}(\Sigma) \in H^2(X;\Z)$, it follows that
$\alpha \cup \beta = 0$.  This implies that $\beta = \alpha \cup
\alpha'$ for some $\alpha' \in H^1(X;\Z)$. But such a $\beta$ satisfies
$\beta\cup \beta = 0$, contradicting the fact that  $\beta \cup \beta
= \Sigma \cdot \Sigma \neq 0$.

To make use of this observation, let $M$ and $M'$ be submanifolds of $X$
dual to $\alpha$. We would like $M$ and $M'$ to be disjoint, so that they
will cobound a submanifold with trivial intersection form.  Consider the
$2$-fold covering space $X_\alpha$ dual to $\alpha$. By
Proposition~\ref{P:toruscover}, we again have a $\Z_{(2)}$ cohomology
torus, and we can iterate, getting a $2^k$ fold covering $X^k\to X$. For
some $k$, there are disjoint lifts of $M$ and $M'$ to $X^k$, and so we
get a spin cobordism from $M$ to $M'$ inside $X^k$, with trivial
intersection form.  Hence the Rohlin invariants of $M$ and $M'$ must
coincide.
\end{proof}


\subsection{Rohlin invariants and the cup product}\label{S:r-cup}
A homology $4$--torus $X$ has an assortment of Rohlin invariants
$\rho(X,\alpha,\sigma)$, associated to different elements $\alpha \in
H^1(X;\Z)$ and different spin structures on $X$. In this subsection,
we show that a certain combination of these Rohlin invariants is
related to the parity of the cup product form. This is, in effect,
half of the main theorem of the paper. The other half will be
established in Section \ref{S:main} by relating the Donaldson
invariants of $X$ to the same cup product.

Denote by $\sigma \to \sigma + x$ the action of $x \in H^1(X;\Z_2)$
on $\sigma \in \Spin(X)$.

\begin{definition}\label{D:rhobar}
Let $X$ be an oriented homology $4$--torus, and let $\alpha \in
H^1(X;\Z)$ be a primitive element such that $\tilde X_\alpha$ has
finitely generated rational homology.  Let $x_1,x_2,x_3 \in
H^1(X;\Z_2)$ be such that $\{\alpha, x_1,x_2,x_3\}$ is a basis for
$H^1(X;\Z_2)$.  Choose a spin structure $\sigma$ on $X$, and define
\begin{equation}\label{E:rhobardef}
\rhobar(X,\alpha) = \sum_{x \in \Span\{x_1,x_2,x_3 \}}
\rho(X,\alpha,  \sigma +x).
\end{equation}
\end{definition}

\begin{lemma}\label{L:welldef}
The invariant $\rhobar(X,\alpha)$ does not depend on the choice of
spin structure $\sigma$, or the specific choice of $\{x_1,x_2,x_3\}$.
\end{lemma}

\begin{proof} The main point in both statements is the following
equality:
$$
\rho(X,\alpha,\sigma) = \rho(X,\alpha, \sigma +\alpha).
$$
This is clear from the definitions, because $\alpha$ restricts
trivially to $M$ and so changing $\sigma$ to $\sigma +\alpha$ does
not affect the spin structure on $M$ (and {\it a fortiori} its Rohlin
invariant).

Given another choice of basis $\{\alpha, x'_1, x'_2, x'_3\}$ for
$H^1 (X;\Z_2)$, the projection along $\alpha$ gives rise to a bijective
correspondence between $\Span\{x_1, x_2, x_3\}$ and $\Span\{x'_1,x'_2,
x'_3\}$. This proves the independence of $\bar\rho(X,\alpha)$ of the
choice of $\{x_1, x_2, x_3\}$. To see that $\bar\rho(X,\alpha)$ does
not depend on $\sigma$, notice that for any two choices, $\sigma$ and
$\sigma'$, of spin structures on $X$,
\begin{align*}
\sum_{x\in\Span\{x_1,x_2,x_3\}}\rho(X,\alpha,\sigma' + x)
&= \sum_{x\in\Span\{x_1,x_2,x_3\}}\rho(X,\alpha,\sigma + (\sigma'
-\sigma) + x) \\
&= \sum_{x'\in\Span\{x'_1,x'_2,x'_3\}}\rho(X,\alpha,\sigma + x'),
\end{align*}
where $\Span\{x'_1,x'_2,x'_3\}$ is the bijective image of $(\sigma'
- \sigma) + \Span\{x_1,x_2,x_3\}$ projected along $\alpha$.
\end{proof}

\begin{remark}
On the other hand, it is not so clear that $\rhobar(X,\alpha)$ is
independent of $\alpha$.  In many cases, this follows from the main
theorem, but we do not know how to deduce this from first principles.
\end{remark}

The relation between the invariant $\bar\rho$ and the cup-product is
established via the following result of Turaev~\cite{turaev:linking},
extending a result of S.~Kaplan~\cite{kaplan:even}.

\begin{theorem}\label{T:kaplancup}
Let $(M^3,\sigma)$ be a spin manifold, and let $x,y,z\in H^1(M;\Z_2)$
be the mod $2$ reductions of integral classes.
Then
\begin{multline}\label{sa}
\rho(M,\sigma) - \rho(M,\sigma + x) - \rho(M,\sigma +y ) -
\rho(M,\sigma + z ) \\ + \rho(M,\sigma + x + y )
+ \rho(M,\sigma + x + z ) + \rho(M,\sigma + y + z ) \\
- \rho(M,\sigma + x + y + z) = (x\cup y\cup z)\,[M]\pmod{2}
\end{multline}
\end{theorem}

\begin{corollary}\label{C:rhocup}
Let $X$ be an oriented homology $4$--torus, and let $\alpha\in H^1(X;\Z)$
be a primitive element such that $\tilde X_\alpha$ is a homology $3$--torus.
Then $\rhobar(X,\alpha) \equiv \det X \pmod 2$.
\end{corollary}

\begin{proof}
This is immediate from Theorem \ref{T:kaplancup} after we notice that all
the Rohlin invariants in (\ref{sa}) take values in $\Z/2\Z$ (by Theorem
\ref{T:rhowelldef}) and that, for any choice of basis
$\{\alpha,x_1,x_2,x_3\}$ in $H^1 (X;\Z_2)$, we have
$(x_1\cup x_2\cup x_3)\,[M] = (\alpha\cup x_1\cup x_2\cup x_3)\,[X] =
\det X\pmod 2$.
\end{proof}


\section{The Donaldson invariants}\label{S:donaldson}

In this section we introduce the $\bar\lambda$--invariants of
$\Z[\Z]$--homology $4$--tori via counting projectively flat connections
on admissible $U(2)$--bundles. The definition is similar to that of the
Donaldson invariants, but with two important differences. First, the
metric perturbations that are usually used in the definition of Donaldson
invariants have no effect here, so we make use of holonomy
perturbations as employed in~\cite{floer:instanton,taubes:casson,
donaldson:orientation,kronheimer:higher-rank} to obtain zero dimensional
perturbed moduli spaces whose points may be counted (with signs). We make
use of the standard cobordism resulting from a homotopy of perturbations
to show that  the invariants $\bar\lambda$ are independent of the
perturbations.  Second, we make essential use in our counting arguments
of an action of $H^1(X;\Z_2)$ on the space of connections, and must take
care that the perturbations we use are equivariant with respect to this
action.  A subtle point here is that we do {\em not} need to use the
stronger claim that the action of $H^1(X;\Z_2)$ extends to the cobordism
constructed from a homotopy of perturbations, compare with Remark
\ref{R:no-equiv}. This stronger statement would be interesting to
establish, but is beyond the methods of the current paper.


\subsection{Admissible bundles}\label{S:admissible}
Let $X$ be a homology $4$--torus and consider a principal $U(2)$--bundle $P$
over $X$ and its associated $SO(3) = PU(2)$--bundle $\bar P$. The
characteristic classes of $P$ and $\bar P$ are related by the formulas
\[
w_2(\bar P) = c_1 (P)\pmod 2\quad\text{and}\quad
p_1(\bar P) = c_1 (P)^2 - 4c_2(P).
\]
Since $H^2(X;\Z)$ is torsion free and $p_1(\bar P) = w_2(\bar P)^2\pmod 4$,
every $SO(3)$--bundle over $X$ arises as $\bar P$ for some $U(2)$--bundle
$P$.

A $U(2)$--bundle $P$ over a homology $4$--torus $X$ is called
\emph{admissible} if $p_1(\bar P) = 0$, $w_2(\bar P)\ne 0$, and, in case
$X$ is even, there exists a $\xi\in H^1(X;\Z_2)$ such that $w_2 (\bar P)
\cup\xi\ne 0$. Note that, for an admissible bundle, $w_2(\bar P)\cup
w_2(\bar P) = p_1(\bar P) = 0 \pmod 4$.  Also note that according to
Corollary \ref{C:cup12}, admissible bundles exist over any
$\Z[\Z]$--homology $4$--torus.

The existence of $\xi \in H^1 (X;\Z_2)$ such that $w_2(\bar P)\cup\xi \ne
0$ is automatic for odd homology tori (by Poincar\'e duality), therefore,
any $U(2)$--bundle $P$ over an odd homology torus such that $p_1(\bar P)=
0$ and $w_2 (\bar P) \ne 0$ is admissible. This is not the case for even
homology tori, as can be seen from the following example.

Consider an even homology $3$--torus, for example $M = \#_3\,(S^1 \times
S^2)$ or the result of $0$--surgery on a band-sum of the Borromean rings
as in~\cite[\S5]{ruberman-strle:tori}.   Note that the cup-product
$\Lambda^2 H^1(M;\Z_2) \to H^2(M;\Z_2)$ is trivial and that $X = M \times
S^1$ is an even $\Z[\Z]$--homology $4$--torus. Choose a nonzero class
$\alpha\in H^1(M;\Z_2)$, and let $\delta$ be the class pulled back from
the $S^1$ factor.  Now let $\bar P$ be the bundle with $p_1(\bar P) = 0$
and $w_2(\bar P) = \alpha \cup \delta$.  No $U(2)$--bundle $P$ lifting
$\bar P$ is admissible because any class $\xi\in H^1(X;\Z_2)$ has trivial
cup product with $w_2 (\bar P)$.


\subsection{$U(2)$ vs $SO(3)$ connections}
Let $X$ be a $\Z[\Z]$--homology $4$--torus and fix a Riemannian metric on $X$.
Every connection on a $U(2)$--bundle $P$ over $X$ induces connections on
$\bar P$ and on the line bundle $\det P$, according to the splitting $\u(2)
= \su(2)\,\oplus\,\u(1)$.  If $A$ is a $\u(2)$--valued 1--form representing
the connection in a local trivialization, this corresponds to the
decomposition
\[
A = \left(A - \frac 1 2\,\tr A\cdot\Id\right) + \frac 1 2\,\tr A\cdot\Id.
\]
(In our discussion, we will denote by $A$ either a connection, or its
local representation as a $1$-form; most calculations will involve
the local representation.)
The induced connection on $\bar P$ is the image of the first summand under
the isomorphism $\ad:\su(2)\to\so(3)$ given by $\ad(\xi)(\eta)=[\xi,\eta]$,
and the induced connection on $\det P$ is $\tr A$.  Conversely, any two
connections on $\bar P$ and $\det P$ determine a unique connection on $P$.

Fix a connection $C$ on $\det P$ and let $\A(P)$ be the space of connections
on $P$ of Sobolev class $L^2_l$ with $l \ge 3$, compatible with $C$. The
assumption $l \ge 3$ ensures that $L^2_l \subset C^0$. The gauge
group $\G(P)$ consisting of unitary automorphisms of $P$ of Sobolev class
$L^2_{l+1}$ having determinant one preserves $C$ and hence acts on
$\A(P)$ with the quotient space $\B(P) = \A(P)/\G(P)$. Let $\A (\bar P)$ be
the affine space of connections on $\bar P$ and $\G(\bar P)$ the full $SO(3)$
gauge group.  Denote $\B(\bar P) = \A(\bar P)/\G(\bar P)$. The natural
projection $\pi:\A(P)\to\A(\bar P)$ commutes with the above gauge group
actions and hence defines a projection
\begin{equation}\label{E:pi}
\pi: \B(P) \to \B(\bar P).
\end{equation}
The space $\B(P)$ admits an $H^1 (X;\Z_2)$ action whose orbit map is the
map (\ref{E:pi}). The rest of this subsection is devoted to a detailed
description of this action, compare with \cite[Section 5.6]{donaldson:floer}.

Let $\O_2: \G(\bar P) \to H^1 (X;\Z_2)$ be the obstruction to lifting
an automorphism $\bar g: \bar P \to \bar P$ to an $SU(2)$ gauge
transformation. Under the standard identification of $H^1 (X; \Z_2)$ as
$\Hom(\pi_1 X,\Z_2)$, the evaluation of $\O_2(\bar g)$ on a loop $\psi$
may be calculated as follows. View $\bar g \in \G(\bar P)$ as a section
of the bundle $\Ad\bar P = \bar P\times_{\Ad} SO(3)$.  The bundle
$\Ad\bar P$ is trivial over $\psi$ so the restriction of $\bar g$
on $\psi$ can be viewed as simply a map to $SO(3)$. Then $\O_2 (\bar g)
(\psi) = \bar g_*(\psi) \in \pi_1(SO(3)) = \Z_2$.
Obstruction theory can be used to show that $\O_2: \G(\bar P) \to
H^1 (X;\Z_2)$ is surjective.

\begin{lemma}
Every gauge transformation $\bar g: \bar P \to \bar P$ lifts to a gauge
transformation $g: P \to P$ of the $U(2)$--bundle $P$.
\end{lemma}

\begin{proof}
Observe that the adjoint representation $U(2)\to SO(3)$ gives rise to a
$U(1)$ bundle $\pi: P \to \bar P$, and view gauge transformations
$g: P \to P$ as sections of the bundle $\Ad  P = P\times_{\Ad} U(2)$.
The result will follow after we calculate the obstruction to lifting a
section $\bar g: X \to \Ad \bar P$ to a section $g: X \to \Ad P$,  and
then show that this obstruction vanishes.

Since the fiber of $\pi: \Ad P \to \Ad \bar P$ is $U(1)$, there is a single
obstruction $o(\bar g)\in H^2(X;\Z)$ to the existence of a lift. By a
standard argument, this obstruction is the pullback, via $\bar g: X \to
\Ad \bar P$, of the obstruction $\O \in H^2(\Ad\bar P;\Z)$ to the existence
of a section of $\pi$.

In order to evaluate $\O$, we apply the Leray--Serre spectral sequence to the
$SO(3)$--bundle $\Ad\bar P \to X$ to calculate the second cohomology of
$\Ad\bar P$. The only relevant non--trivial differentials are those at the
$E_3$ level, and these would be $d_3: E_3^{0,2} \to E_3^{3,0}$.  But the
former group is $H^2(SO(3);\Z) = \Z_2$ and the latter is $H^3(X;\Z)$,
which is torsion free. Thus we have an exact sequence
\[
\begin{CD}
0 @>>> H^2(X;\Z) @>p^*>> H^2(\Ad\bar P;\Z) @>i^*>> H^2(SO(3);\Z) @>>> 0,
\end{CD}
\]
where $i: SO(3) \to \Ad\bar P$ is the inclusion of $SO(3)$ as the fiber.
This sequence is split because there is a section of the bundle $\Ad\bar P
\to X$, which is the identity gauge transformation $\Id_{\bar P}$.
Thus $H^2(\Ad\bar P;\Z) \cong H^2(X;\Z) \oplus H^2(SO(3);\Z)$.   In
particular, since the cohomology of $X$ is torsion free, there is a
unique $2$--torsion class $t \in H^2 (\Ad\bar P;\Z)$, and this class
restricts to $H^2(SO(3);\Z)$ as the generator.

We claim that $\O = t$. Since the above exact sequence splits, it suffices
to show that $\Id_{\bar P}^*(\O) = 0$ but that $\O$ is nontrivial. Note that
the identity gauge transformation $\Id_{\bar P}$ lifts to the identity
transformation of $P$, hence $\Id_{\bar P}^*(\O)$, which is the obstruction
to this lifting, must vanish. On the other hand, we have that $i^*\O \neq 0
\in H^2(SO(3);\Z)$ because this pullback represents the (complete)
obstruction to a section of the adjoint representation $U(2) \to SO(3)$.
But there certainly is no such section.

Finally, $o(\bar g) = \bar g^*(\O) $ is a $2$--torsion element of $H^2(X;\Z)$
(since $\O$ is $2$--torsion). Since $H^2(X;\Z)$ is torsion free, $o(\bar g) =
0$ and $\bar g$ admits a lift.
\end{proof}

Given $\chi \in H^1 (X; \Z_2)$, choose $\bar g \in \G(\bar P)$ so that
$\O_2(\bar g) = \chi$. Any gauge transformation $g: P \to P$ lifting
$\bar g$ will be said to \emph{realize} $\chi$.

\begin{lemma}\label{L:lift}
Let $\chi \in H^1 (X;\Z_2)$ be realized by $g: P \to P$ and let $h = \det g:
X \to S^1$.  Then $h$ induces a homomorphism $h_* \in \Hom\, (\pi_1 X, \Z) =
H^1 (X; \Z)$ whose modulo 2 reduction is $\chi$.
\end{lemma}

\begin{proof}
The result follows by a direct calculation with loops $\psi \in \pi_1 X$
using the observations that $\det: U(2) \to  U(1)$ induces an isomorphism
of the fundamental groups $\pi_1\,U(2) = \pi_1\,U(1) = \Z$ and that the
homomorphism  $\pi_1\,U(2) \to \pi_1\,SO(3)$ induced by the adjoint
representation $U(2) \to SO(3)$ is the mod 2 reduction $\Z \to \Z_2$.
\end{proof}

Given $\chi \in H^1 (X;\Z_2)$, realize it by a gauge transformation $g: P
\to P$ as above and define the action of $\chi$ on $\B(P)$ by the formula
\begin{equation}\label{E:chi}
\chi[A] = [\,g^* A - 1/2\,h\,dh^{-1}\Id\,]
\end{equation}
where $h = \det g$ and, in local trivializations, $g^* A = g\,dg^{-1} +
g A g^{-1}$.  Since $\tr\,(g^* A - 1/2\,h\,dh^{-1}\Id) = \tr A +
\tr\,(g\,dg^{-1}) - h\,dh^{-1} = \tr A$, the determinant connection is
preserved and hence $\chi[A]$ belongs to $\B(P)$.

Any two gauge transformations $g_1, g_2: P \to P$  lifting $\bar g: \bar P
\to \bar P$ differ by an automorphism $\gamma = c \Id$,  where $c\in U(1)$,
so that $g_2 = \gamma g_1$.  A straightforward calculation then shows that
\begin{alignat}{1}
g_2^* A - 1/2\,h_2\,dh_2^{-1}\Id
&= (\gamma g_1)^* A - 1/2\;c^2 h_1\,d(c^{-2} h_1^{-1})\,\Id \notag \\
&= g_1^* A + \gamma d\gamma^{-1} - 1/2\,c^2 h_1 (2c^{-1}dc^{-1}h_1^{-1} +
c^{-2} dh_1^{-1})\Id \notag \\
&= g_1^* A + \gamma\,d\gamma^{-1} - c\,dc^{-1}\,\Id - 1/2\,h_1\,dh_1^{-1}
\Id \notag \\
&= g_1^* A - 1/2\, h_1\,dh_1^{-1}\Id,\label{E:calc}
\end{alignat}
which implies that (\ref{E:chi}) is independent of the choice of a lift of
$\bar g: \bar P \to \bar P$.  Furthermore, any two choices of $\bar g$
realizing the same $\chi$ differ by an $SU(2)$ gauge transformation, and so
do their respective lifts $g_1, g_2: P \to P$. But then the connections
$g_2^* A - 1/2\,h_2\,dh_2^{-1}\Id$ and $g_1^* A - 1/2\,h_1\,dh_1^{-1}\Id$
differ by an $SU(2)$ gauge transformation and thus define the same element
in $\B(P)$. Therefore, (\ref{E:chi}) is independent of the choices made in
its definition.

Given $\chi_1$ and $\chi_2 \in H^1 (X;\Z_2)$ realized by their respective
lifts $g_1$ and $g_2: P\to P$, the product $\chi_1\chi_2$ can be realized
by the composition $g_1 g_2: P \to P$. In particular,
\begin{multline}\notag
\chi_1(\chi_2 [A])
= [\,g_1^* (g_2^* A - 1/2\,h_2\,dh_2^{-1}\Id) - 1/2\,h_1\,dh_1^{-1}\Id\,] \\
= [\,(g_1 g_2)^* - 1/2\,h_1 h_2\,d(h_1 h_2)^{-1}\Id\,] = (\chi_1\chi_2)[A]
\end{multline}
since $h_1 = \det g_1$ and $h_2 = \det g_2$ commute. Moreover, $\chi = 0$
can be realized by $\Id: P \to P$ and hence induces an identity map on
$\B(P)$. This completes the definition of the $H^1 (X;\Z_2)$ action on
$\B(P)$.


\subsection{Projectively ASD connections}\label{S:proj-asd}
A connection $A$ on an admissible bundle $P$ is called \emph{projectively
ASD} if the connection $\bar A$ it induces on the bundle $\bar P$ is ASD.
The latter means that $F_+ (\bar A) = 0$ where, as usual, $F_+ (\bar A) =
1/2\cdot(F(\bar A) + *\,F(\bar A))$ is the self dual part of the curvature.
Equivalently, a connection $A$ on an admissible bundle $P$ is projectively
ASD if $F_+(A)$ is central. Once a connection $C$ on $\det P$ is fixed,
the space of all projectively ASD connections on $P$ compatible with $C$,
modulo gauge group $\G(P)$, will be denoted by $\M(P)$. It is called the
\emph{moduli space of projectively ASD connections}; different choices of
$C$ give equivalent moduli spaces.  The moduli space $\M(P)$ (and its
perturbed version, $\M_\sigma(P)$, defined below) depends on the
choice of Riemannian metric on $X$.  When we want to indicate this
dependence, we will use the notation $\M_g(P)$ (or
$\M_{g,\sigma}(P)$).

The ASD equation is equivariant with respect to the $H^1 (X;\Z_2)$ action
in the following sense. The group $\G(P)$ acts on $\Omega^2_+ (X; \ad\bar
P)$ completed in Sobolev $L^2_{l-1}$ norm by pull back of differential
forms, an element $g\in \G(P)$ sending $\omega$ to $(\ad\bar g)^* \omega$.
Let
\begin{equation}\label{E:bndl}
\E(P) = \A(P)\, \times_{\,\G(P)}\, \Omega^2_+ (X;\ad\bar P)
\end{equation}
be the vector bundle over $\B(P)$ associated with the principal $\G(P)$
bundle $\A(P) \to \B(P)$. One can easily check that the formula $F_+[A]
= [A, F_+ (\bar A)]$ defines a smooth section $F_+: \B(P) \to \E(P)$ of
the above bundle.

\begin{proposition}\label{P:action}
There is a natural action of $H^1(X;\Z_2)$ on the bundle $\E(P)\allowbreak
\to \B(P)$ extending the action on $\B(P)$ such that the section $F_+$ is
equivariant.
\end{proposition}

\begin{proof}
Given a $\chi \in H^1 (X; \Z_2)$ and a gauge transformation $g: P \to P$
realizing it, define $\chi: \E(P)\to \E(P)$ by the formula $\chi[A,\omega]
=  [g^* A - 1/2\,h\,dh^{-1}\Id,\allowbreak (\ad\bar g)^*\omega]$ where $h
= \det g$. That this definition is independent of the arbitrary choices
follows by essentially the same argument as well--definedness of the action
(\ref{E:chi}) on $\B(P)$. The equivariance of $F_+$ is immediate from the
observation that $h\,dh^{-1}$ is a closed 1--form and hence $F(g^* A - 1/2
\,h\,dh^{-1}\Id) = (\ad g)^* F(A)$.
\end{proof}

\begin{corollary}\label{C:chi}
The moduli space $\M(P)$ is the zero set of the equivariant section $F_+:
\B(P) \to \E(P)$ and as such it is acted upon by $H^1 (X; \Z_2)$.  The
quotient space of this action is the moduli space $\M(\bar P)$ of ASD
connections on $\bar P$ modulo gauge group $\G(\bar P)$, compare with
(\ref{E:pi}).
\end{corollary}

Let $A \in \A(P)$ then the map $\A(P) \to \Omega^1 (X; \ad\bar P)$ sending
$B$  to  $\bar B - \bar A$ establishes an isomorphism between  $T_A \A(P)$
and $\Omega^1(X;\ad\bar P)$ completed in Sobolev $L^2_l$ norm. Under this
isomorphism, the slice through $A$ of the $\G(P)$ action on $\A(P)$ is
isomorphic to $\ker d^*_{\bar A} \subset \Omega^1(X;\ad\bar P)$. The map
$F_+: \B(P) \to \E(P)$ linearizes at $[A]$ to $d^+_{\bar A}:
\ker d^*_{\bar A} \to \Omega^2_+ (X;\ad\bar P)$ so the local structure of
$\M(P)$ near $[A]$ is described by the deformation complex
\begin{equation}\label{D:full}
\begin{CD}
\Omega^0 (X;\ad\bar P) @> d_{\bar A} >> \Omega^1 (X;\ad\bar P)
@> d^+_{\bar A} >> \Omega^2_+ (X;\ad\bar P)
\end{CD}
\end{equation}
with cohomology $H^0 (X;\ad A)$, $H^1 (X;\ad A)$, and $H^2_+ (X;\ad A)$.
Note for future use that the deformation complexes at $[A]$ and $\chi[A]
\in \M(P)$ are isomorphic to each other via $(\ad\bar g)^*$ where $g$
realizes $\chi$\,:
\begin{equation}\label{D:iso}
\begin{CD}
\Omega^0 (X;\ad\bar P) @> d_{\bar g^*\bar A} >> \Omega^1 (X;\ad\bar P)
@> d^+_{\bar g^*\bar A} >> \Omega^2_+ (X;\ad\bar P) \\
@VV (\ad\bar g)^* V @VV (\ad\bar g)^*V @VV (\ad\bar g)^*V \\
\Omega^0 (X;\ad\bar P) @> d_{\bar A} >> \Omega^1 (X;\ad\bar P)
@> d^+_{\bar A} >> \Omega^2_+ (X;\ad\bar P).
\end{CD}
\end{equation}
The above allows for a computation of the formal dimensions of both $\M(P)$
and $\M(\bar P)$ using the index theorem\,: since $X$ is a $\Z[\Z]$--homology
$4$--torus, we have $b_2^+(X) = 3$ and therefore
\[
\dim \M(P) = \dim \M(\bar P) = -2 p_1 (\bar P) - 3(1 - b_1 + b_2^+)(X) = 0.
\]
A straightforward application of the Chern--Weil theory to the bundle
$\bar P$ with $p_1(\bar P) = 0$ shows that projectively ASD connections
on $P$ are in fact \emph{projectively flat}, which means that they induce
flat connections on $\bar P$. Therefore, the moduli spaces $\M(P)$ and
$\M(\bar P)$ can be viewed as moduli spaces of projectively flat
(respectively, flat) connections. The holonomy correspondence then
identifies $\M(\bar P)$ with a compact subset of the $SO(3)$--character
variety of $\pi_1 X$. Since $\M(P)$ is a finite sheeted covering of
$\M(\bar P)$, see Corollary \ref{C:chi}, it is also compact. Proposition
\ref{P:hol} below will provide an explicit description of $\M(P)$ in terms
of (projective) $SU(2)$--representations of $\pi_1 X$.

We say that $\M(P)$ is \emph{non--degenerate} at $[A] \in \M(P)$ if $H^2_+
(X;\ad A) = 0$.  According to (\ref{D:iso}), $\M(P)$ is non--degenerate at
$\chi[A]$ if and only if it is non--degenerate at $[A]$. We say that
$\M(P)$ is \emph{non--degenerate} if it is non--degenerate at all $[A] \in
\M(P)$. The above discussion implies that, if $\M(P)$ is non--degenerate
then it is a smooth compact manifold of dimension zero. It is canonically
oriented once a homology orientation is fixed, see
\cite[Section 3]{donaldson:orientation}, and the action of $H^1(X;\Z_2)$
on $\M(P)$ is then orientation preserving, see
\cite[Corollary 3.27]{donaldson:orientation}.

If $\M(P)$ fails to be non--degenerate, the same holds true after $\M(P)$
is properly perturbed.


\subsection{Admissible perturbations}\label{S:adm-pert}
Let $S^1$ be the interval $[0,1]$ with identified ends, and consider an
embedding $\psi: S^1\times D^3 \to X$ so that, for each point $z \in D^3$,
we have a loop $\psi_z = \psi (S^1 \times\{z\})$. The loop $\psi_0$ will
also be called $\psi$. We will say that $\psi$ is \emph{mod-2 trivial} if
$0 = [\psi] \in H_1 (X;\Z_2)$.

Given $A \in \A(P)$, denote by $\hol_A (\psi)$ the function on $\psi(S^1
\times D^3)$ whose value at $\psi(s,z)$ equals the $SU(2)$ holonomy of
$A$ around the loop $\psi_z$ starting at the point $\psi(s,z)$. Recall
from \cite[Section 3.2]{ruberman-saveliev:casson} that, once a connection
$C$ on $\det P$ is fixed, a choice of square roots of $\hol_C$ on a set
of representative loops makes $\hol_A (\psi)$ into a well defined $SU(2)$
valued function, with different choices leading to equivalent theories.
Let $\Pi: SU(2)\to \su(2)$ be the projection given by
\[
\Pi (u) = u - \frac 1 2\,\tr(u)\cdot\Id.
\]
It is equivariant with respect to the adjoint action of $SU(2)$ on both
$SU(2)$ and $\su(2)$,  therefore, $\Pi \hol_A (\psi)$ is a well defined
section of $\ad\bar P$ over $\psi (S^1\times D^3)$ after we identify
$\su(2)$ with $\so(3)$. For any $\nu \in \Omega^2_+ (X)$ supported in
$\psi(S^1\times D^3)$, the formula
\begin{equation}\label{E:sigma}
\sigma(\nu,\psi)[A] = [\,A, \Pi \hol_A(\psi)\,\otimes\,\nu\,]
\end{equation}
defines a section $\sigma(\nu,\psi): \B(P)\to \E(P)$ of the bundle
(\ref{E:bndl}). That $\sigma(\nu,\psi)$ is well defined follows from the
fact that $\Pi\hol_{g^* A}(\psi) = \ad\bar g \cdot \Pi\hol_A (\psi)$ for
any $g \in \G(P)$.

\begin{proposition}\label{P:equiv}
If $\psi$ is mod-2 trivial then the section (\ref{E:sigma}) is equivariant
with respect to the $H^1 (X;\Z_2)$ action on the bundle $\E(P) \to \B(P)$,
compare with Proposition \ref{P:action}.
\end{proposition}

\begin{proof}
Given $\chi \in H^1(X; \Z_2)$, realize it by a gauge transformation $g: P
\to P$. Let $h = \det g: X \to S^1$ and consider the flat connection
$1/2\,h\,dh^{-1}$ on a trivial line bundle $L \to X$. We wish first to
compute the holonomy of this connection.

For an arbitrary loop $\psi: S^1 \to X$, the restriction of $h$ to $\psi$
is a map $S^1 \to S^1$ of the form $t \to \exp(i f(t))$ where $f(2\pi) =
f(0) + 2\pi k$ and $k = h_*(\psi) \in \Z$. Then $1/2\,h\,dh^{-1} = i/2\,
f'(t)\,dt$, and the holonomy of $1/2\,h\,dh^{-1}$ around $\psi$ equals
\smallskip
\[
\exp\left(\frac i 2\,\int_0^{2\pi} f'(t)\,dt\right) =
\exp\left(\frac i 2\,(f(2\pi) - f(0))\right) = \exp(i \pi k) = \chi(\psi),
\]
where $\chi$ is viewed as a homomorphism from $\pi_1 X$ to $\Z_2 =
\{\,\pm 1\,\}$, compare with Lemma \ref{L:lifts}.

Let $\psi$ be a mod-2 trivial loop. For any $A \in \A(P)$, the connection
$g^*A - 1/2\,h\,dh^{-1}\Id =  g^*\,(A - 1/2\,h\,dh^{-1}\Id)$ can be viewed
as the pull back via $g: P \to P$ of the connection on $P\,\otimes\,L = P$
induced by $A$ and $-1/2\,h\,dh^{-1}$. Then
\begin{alignat*}{1}
\Pi\hol_{g^* A - 1/2\,h\,dh^{-1}\Id}(\psi)
&= \ad\bar g\cdot \Pi(\hol_A(\psi)\cdot\hol_{-1/2\,h\,dh^{-1}}(\psi))  \\
&= \ad\bar g\cdot \Pi\,(\hol_A(\psi)\cdot\chi(\psi))
= \ad\bar g\cdot \Pi\hol_A(\psi),
\end{alignat*}
after we notice that $\chi(\psi) = 1$ because $\chi:\pi_1 X \to \Z_2$
factors through $H_1(X;\Z_2)$ and $\psi$ is mod-2 trivial. Therefore,
\[
\Pi\hol_{g^* A - 1/2\,h\,dh^{-1}\Id}(\psi)\otimes\nu = (\ad\bar g)^*
(\Pi\hol_A(\psi)\otimes\nu),
\]
which implies equivariance of the section $\sigma(\nu,\psi)$.
\end{proof}

An \emph{admissible perturbation} is a smooth section $\sigma$ of the
bundle $\E(P) \to \B(P)$ of the form
\begin{equation}\label{E:sigma'}
\sigma = \sum_{k=1}^N\;\ep_k\cdot\sigma(\nu_k,\psi_k), \quad \ep_k \in
\mathbb R,
\end{equation}
where $\sigma(\nu_k,\psi_k)$ are sections (\ref{E:sigma}) built from
a family of loops $\psi_k: S^1 \times D^3 \to X$ with disjoint
images, and $\nu_k \in \Omega^2_+ (X)$ are self--dual $2$--forms, each
supported in its respective $\psi_k (S^1\times D^3)$.  An admissible
perturbation is \emph{equivariant} if it is equivariant with respect
to the action of $H^1 (X;\Z_2)$ on the bundle $\E(P) \to \B(P)$.
According to Proposition \ref{P:equiv}, all admissible perturbations
built from mod-2 trivial loops $\psi_k$ are equivariant.

Let $\sigma: \B(P) \to \E(P)$ be an admissible perturbation and define
$\zeta_{\sigma}: \B(P) \to \E(P)$ by the formula $\zeta_{\sigma} = F_+
+ \sigma$. We call a connection $A \in \A(P)$ \emph{perturbed
projectively ASD} if $\zeta_{\sigma}\,[A] = 0$, and denote by
$\M_{\sigma}(P) = \zeta^{-1}_{\sigma}(0)$ the moduli space of
perturbed projectively ASD connections. If $\sigma = 0$ then
$\M_{\sigma}(P)$ is the moduli space $\M(P)$ of projectively ASD
(equivalently, projectively flat) connections.

Any admissible perturbation $\sigma: \B(P) \to \E(P)$ factors through
the embedding $\E'(P) \to \E(P)$ where $\E'(P)$ is the bundle
(\ref{E:bndl}) with fiber $\Omega^2_+ (X;\ad\bar P)$ completed in
Sobolev $L^2_l$ norm, see \cite[Proposition 7]{kronheimer:higher-rank}.
Since the inclusion $L^2_l\to L^2_{l-1}$ is compact, the derivatives
of $\sigma$ are compact operators and hence the perturbed section
$\zeta_{\sigma}$ is Fredholm of index zero.

\begin{proposition}\label{P:compact}
For every admissible perturbation $\sigma$, the perturbed moduli space
$\M_{\sigma}(P)$ is compact.
\end{proposition}

\begin{proof}
The key to proving compactness is the observation that the forms
$\Pi\hol_A(\psi)\,\otimes\,\nu$ admit a uniform $L^\infty$ bound
independent of $A$. Given a sequence $[A_n] \in \M_{\sigma} (P)$, the
perturbed projectively ASD equations then imply that the $F_+ (\bar A_n)$
are uniformly bounded in $L^\infty$, and $F (\bar A_n)$ are uniformly
bounded in $L^2$ because of the Chern--Weil formula. Using Uhlenbeck's
compactness theorem \cite{uhlenbeck:lp-bounds} in the absence of bubbling, we
conclude that (after passing to a subsequence and gauge equivalent
connections, if necessary) the sequence $A_n$ converges in $L^p_1$ for
all $p \ge 2$ to a connection $A$ such that $[A] \in \M_{\sigma}(P)$,
compare with \cite[Proposition 11]{kronheimer:higher-rank}. Finally,
bootstrapping leads to the $C^\infty$ convergence.
\end{proof}

The local structure of $\M_{\sigma}(P)$ near a point $[A]$ is described
by the deformation complex
\begin{equation}\label{D:pert}
\begin{CD}
\Omega^0 (X,\ad\bar P) @> d_{\bar A} >> \Omega^1 (X,\ad\bar P)
@> d^+_{\bar A} + D_A \sigma >> \Omega^2_+ (X,\ad\bar P).
\end{CD}
\end{equation}
Here, assuming that $\sigma$ is given by the formula (\ref{E:sigma'}),
we have $D_A \sigma = \sum\;\ep_k\cdot D_A \sigma(\nu_k, \psi_k)$
where $D_A \sigma(\nu_k,\psi_k)$ is the differential at $A$ of the
function $\A(P)\to \Omega^2_+ (X; \ad\bar P)$ mapping $A$ to
$\Pi\hol_A(\psi_k) \otimes \nu_k$.

We call $\M_{\sigma}(P)$ \emph{non--degenerate} at $[A]\in \M_{\sigma}
(P)$ if the second cohomology of the complex (\ref{D:pert}) vanishes,
that is, $\coker\,(d_{\bar A}^+ + D_A \sigma) = 0$. We say that
$\M_{\sigma}(P)$ is \emph{non--degenerate} if it is non--degenerate at
all $[A] \in \M_{\sigma}(P)$. A non--degenerate $\M_{\sigma}(P)$ is a
compact zero--dimensional manifold,  canonically oriented by a choice
of homological orientation. If $\sigma$ is equivariant then,
according to Proposition \ref{P:equiv}, $\M_{\sigma} (P)$ is acted
upon by $H^1(X; \Z_2)$, and this action is orientation preserving.

Our immediate goal will be to show that $\M_{\sigma}(P)$ is
non--degenerate for a generic (not necessarily equivariant) admissible
perturbation $\sigma$. This will allow us to define the Donaldson
invariant in Section \ref{S:l-def};  later in the paper, we will prove
a more refined result that $\M_{\sigma}(P)$ can be made non--degenerate
using a generic \emph{equivariant} admissible perturbation $\sigma$.

Recall that, after making certain choices as in \cite[Section 3.2]
{ruberman-saveliev:casson}, the holonomy of a connection  $A \in \A(P)$
can be viewed as a map $\hol_A: \Omega (X,x) \to SU(2)$ on the monoid
of loops based at $x \in X$ such that $\hol_A (\psi_1 \psi_2) = \pm\,
\hol_A(\psi_1)\,\hol_A(\psi_2)$. We will call $A$ \emph{reducible} if
the image of $\hol_A$ is contained in a copy of $U(1)\subset SU(2)$,
and \emph{irreducible} otherwise.  Among reducible connections,  those
with the image of $\hol_A$ in the center of $SU(2)$ will be called
\emph{central}.

\begin{lemma}\label{L:don25}
If $A$ is a projectively flat irreducible connection on $P$ then there
are finitely many loops $\psi_k: S^1 \times D^3 \to X$ with
disjoint images and self--dual $2$--forms $\nu_k$, each supported in a
small ball inside $\psi_k(S^1\times D^3)$, such that the sections
$\Pi\hol_A(\psi_k)\,\otimes\,\nu_k \in \Omega^2_+(X;\ad\bar P)$
generate the vector space $H^2_+(X;\ad A)$.
\end{lemma}

\begin{proof}
This is essentially Lemma 2.5 from Donaldson \cite{donaldson:orientation}.
\end{proof}

The moduli space $\M(P)$ consists of the gauge equivalence classes of
projectively flat connections all of which are irreducible because of
the assumption that $w_2 (\bar P) \ne 0$. Since $\M(P)$ is compact,
we can suppose that the loops $\psi_k$ and the forms $\nu_k$ can
be chosen so that Lemma \ref{L:don25} holds for all the points $[A]\in
\M(P)$ simultaneously. We fix a choice of $\psi_k$ and $\nu_k$,  and
define $\sigma$ by the
formula (\ref{E:sigma'}).

\begin{proposition}\label{P:gen}
There exists a real number $r > 0$ such that,  with the choice of
$\psi_k$ and $\nu_k$ as above,  the moduli space $\M_{\sigma}(P)$
is non--degenerate for a generic $\bar\ep = (\ep_1,\ldots,\ep_N)$
inside the ball $\|\bar\ep\| < r$ in $\mathbb R^N$.
\end{proposition}

\begin{proof}
Let us pull the bundle $\E(P)$ back to a bundle $\pi^* \E(P) \to
\B(P) \times \mathbb R^N$  via the projection $\pi: \B(P) \times
\mathbb R^N \to \B(P)$, and consider the ``universal section''
\begin{equation}\label{E:psi}
\Psi([A],\bar\ep) = [\,A,\bar\ep, F_+(\bar A) + \sum\;\ep_k \cdot
\Pi\hol_A (\psi_k)\,\otimes\,\nu_k\,].
\end{equation}
Let $[A] \in \M(P)$ and choose $a \in T_{[A]} \A(P) = \ker
d^*_{\bar A}$,\; $\bar\eta = (\eta_1,\ldots,\eta_N) \in
\mathbb R^N$, and a small $t \in \mathbb R$.  Then, up to higher
order terms in $t$, we have
\begin{alignat*}{1}
\Psi([A + ta],t\bar\eta) &= [\,A + ta, t\bar\eta, F_+(\bar A + ta)
+ \sum\; t\eta_k\cdot\Pi\hol_{A + ta}(\psi_k)\otimes\nu_k\,] \\
&= \Psi([A],0) + t\,[\,a,\,\bar\eta,\,d^+_{\bar A} a + \sum\;\eta_k
\cdot \Pi\hol_A(\psi_k)\otimes\nu_k\,] + \ldots
\end{alignat*}
Since sections $\Pi\hol_A(\psi_k) \otimes \nu_k$ generate $H^2_+
(X;\ad A) = \coker d^+_{\bar A}$, we conclude that $\Psi$ is
transversal at $\M(P) \times \{\,0\}$ to the zero section of
$\pi^* \E(P) \to \B(P) \times \mathbb R^N$. Therefore, there
exist an open neighborhood $\U$ of $\M(P)$ in $\B(P)$ and a ball
$\V \subset \mathbb R^N$ centered at the origin such that the
intersection of $\Psi^{-1} (0)$ with the product $\U \times \V
\subset \B(P)\times \mathbb R^N$ is a smooth manifold of dimension
$N$.

Let us consider the projection $p: \Psi^{-1}(0)\cap (\U \times \V)
\to \V$.  For every $\bar\ep \in \V$ we have $p^{-1}(\bar\ep) =
\M_{\sigma}(P)\, \cap\, \U$ where $\sigma$ is given by the formula
(\ref{E:sigma'}). By Uhlenbeck's compactness theorem, if $\bar\ep$
is small enough, then the entire moduli space $\M_{\sigma}(P)$ is
contained in $\U$\,: otherwise, there would exist a sequence $A_n$
of perturbed projectively ASD connections with $\|F (\bar A_n)\|_{L^2}
\to 0$ and $\|F_+(\bar A_n)\|_{L^p} \to 0$ but with no subsequence
converging in $L^p_1$ to a flat connection. Thus there exists an
$r > 0$ such that $p^{-1} (\bar\ep) = \M_{\sigma}(P)$ as long as
$\|\bar\ep\| < r$.  The Sard--Smale Theorem now implies that the
projection $p$ is a submersion for an open dense subset of the ball
$\|\bar\ep\| < r$, that is, that $\coker(d^+_{\bar A} + D_A\sigma)
= 0$ for all $[A]\in \M_{\sigma}(P)$.
\end{proof}


\subsection{Definition of the invariants}\label{S:l-def}
Let $P$ be an admissible $U(2)$--bundle over a $\Z[\Z]$--homology $4$--torus
$X$. After perturbation, if necessary, the moduli space $\M_{\sigma}(P)$
consists of finitely many points,  and it is canonically oriented up to a
choice of homology orientation. We define the invariant $\bar\lambda(X,P)$
as one quarter times signed count of points in $\M_{\sigma}(P)$,
\[
\bar\lambda (X,P) = \frac 1 4\,\#\M_{\sigma}(P).
\]
In what follows, we will only be interested in the modulo 2 reduction
of $\bar\lambda (X,P)$, therefore, a particular choice of homology
orientation, which is implicit in the above definition, will not matter.

The invariant $\bar\lambda(X,P)$ is well defined because, for any generic
path of Riemannian metrics $g(t)$ and small perturbations
$\sigma(t)$, $t \in [0,1]$, the moduli spaces
$\M_{g(0),\,\sigma(0)}(P)$ and $\M_{g(1),\,\sigma(1)}(P)$ are cobordant
via the compact
oriented cobordism
\[
\overline{\M}(P) = \bigcup_{t\in 
[0,1]}\;\M_{g(t),\,\sigma(t)}(P)\times\{\,t\,\}.
\]
It is not hard to show that $\bar\lambda(X,P)$ is
independent of the choice of determinant connection $C$.

\begin{remark}\label{R:no-equiv}
It should be noted that the definition of $\bar\lambda (X,P)$ does not
require that $\sigma$ be equivariant. On the other hand, Sections
\ref{S:eight} and \ref{S:sixteen} provide us with an equivariant
admissible perturbation $\sigma$ such that $\M_{\sigma}(P)$ is
non--degenerate \emph{and} is acted upon by $H^1(X;\Z_2)$. Since this
action is orientation preserving, all points in the orbit of a
(perturbed) projectively flat connection are counted in $\bar\lambda(X,P)$
with the same sign. It would be tempting to define the invariant
$\bar\lambda(X,P)$ by counting points in $\M(\bar P)$ instead of those in
$\M(P)$. However, such an approach comes against two problems\,: one is
that the points of $\M(\bar P)$ would need to be counted with weights
given by the order of the orbits of their respective lifts to $\M(P)$. A
more serious problem would be proving the well definedness of such an
invariant\,: one would have to deal with a generic path of
\emph{equivariant} perturbations, which might not be an easy task.
\end{remark}

\begin{remark}
The class of admissible perturbations used in Sections \ref{S:eight} and
\ref{S:sixteen} is larger than one defined in Section \ref{S:adm-pert}.
However, the above argument showing well--definedness of $\bar\lambda(X,P)$
remains valid for these more general perturbations given their compactness
properties, as described in Section \ref{S:tau}.
\end{remark}


\section{Projective representations}\label{S:proj}

The holonomy map embeds the moduli space $\M(\bar P)$ of flat connections on
$\bar P$ into the $SO(3)$--representation variety of $\pi_1 (X)$. An attempt
to extend this correspondence to the moduli space $\M(P)$ of projectively
flat connections leads to the concept of a projective representation. We
briefly describe the construction and refer the reader
to~\cite{ruberman-saveliev:casson} for a detailed treatment.


\subsection{Algebraic background}\label{S:algebra}

Let $G$ be a finitely presented group and view $\Z_2 = \{\pm 1\}$ as the
center of $SU(2)$. A map $\rho: G\to SU(2)$ is called a \emph{projective
representation} if $\rho(gh)\rho(h)^{-1}\rho(g)^{-1}$ belongs to $\Z_2$ for
all $g, h\in G$. Given a projective representation $\rho$, the function $c:
G\times G\to\Z_2$ defined as $c(g,h) = \rho(gh)\rho(h)^{-1}\rho(g)^{-1}$ is
a $2$--cocycle, that is, $c(gh,k)c(g,h) = c(g,hk)c(h,k)$. We will refer to
$c$ as the \emph{cocycle associated with $\rho$}.

Let us fix a cocycle $c: G\times G\to\Z_2$ and denote by $\P_c (G;SU(2))$
the set of conjugacy classes of projective representations $\rho: G
\to SU(2)$ whose associated cocycle is $c$. Up to an isomorphism, this
set only depends on the cohomology class $[c]\in H^2 (G;\Z_2)$.

The group $H^1 (G;\Z_2)$ acts on $\P_c (G; SU(2))$ by the formula $\rho
\mapsto \rho^{\chi}$ where $\rho^{\chi}(g) = \chi(g)\rho(g)$. The quotient
space of this action is the space $\R_{[c]}(G; SO(3))$ consisting of
(the conjugacy classes of) $SO(3)$--representations whose second
Stiefel--Whitney class equals $[c]$. The projection map $\pi:
\P_c (G; SU(2)) \to \R_{[c]}(G; SO(3))$ can be given explicitly by
composing $\rho: G \to SU(2)$ with the adjoint representation $\ad: SU(2)
\to SO(3)$.

A projective representation $\rho: G\to SU(2)$ is called \emph{irreducible}
if the centralizer of its image equals the center of $SU(2)$. All projective
representations whose associated $2$--cocycle is not homologous to zero are
irreducible, assuming that $H_1 (G;\Z)$ has no $2$--torsion.


\subsection{The holonomy correspondence}\label{S:hol}

In this section we establish a correspondence between projectively flat
connections over a $\Z[\Z]$--homology $4$--torus $X$ and projective
representations of its fundamental group. We first deal with the
discrepancy arising from the fact that $H^2(\pi_1 X;\Z_2)$ and
$H^2(X;\Z_2)$ need not be isomorphic.

\begin{lemma}\label{L:empty}
Let $X$ be a $\Z[\Z]$--homology torus and $P$ an admissible $U(2)$--bundle
over $X$. If $w_2 (\bar P)$ does not belong to the image of $\iota:
H^2(\pi_1 X;\Z_2) \to H^2 (X;\Z_2)$, see (\ref{E:iota}), then the moduli
space $\M(\bar P)$ is empty.
\end{lemma}

\begin{proof}
The Hopf exact sequence $\pi_2 X \to H_2(X;\Z) \to H_2(\pi_1 X;\Z) \to 0$,
see \cite{brown:cohomology}, implies that, if $w_2(\bar P)$ does not 
belong to the
image of $\iota$, it evaluates non-trivially on a $2$--sphere in $X$. Such a
bundle $\bar P$ cannot support any flat connections, for a flat connection
on $\bar P$ would pull back to a flat connection on the $2$--sphere, whose
holonomy would trivialize the bundle.
\end{proof}

\begin{corollary}\label{C:iota}
The main theorem holds for all $\bar\lambda(X,P)$ such that $w_2(\bar P)$
is not in the image of $\iota: H^2(\pi_1 X;\Z_2) \to H^2(X;\Z_2)$.
\end{corollary}

\begin{proof}
Let $P$ be an admissible $U(2)$--bundle  such that $w_2(\bar P)$
is not in the image of $\iota$. Then, according to Lemma
\ref{L:empty}, the moduli space
$\M(\bar P)$ is empty. On the other hand, this situation is only possible
if $X$ is an even homology torus, see Corollary \ref{C:odd}.
\end{proof}

           From now on, we will concentrate on admissible $U(2)$--bundles $P$
such that $w_2 (\bar P)$ is in the image of $H^2(\pi_1 X;\Z_2)$,  and will
identify $H^2(\pi_1 X;\Z_2)$ with its (monomorphic) image in $H^2(X;\Z_2)$.

It is a well known fact that the holonomy defines a bijection $\bar\phi:
\M(\bar P) \to \R_w (X;SO(3))$ where $w = w_2(\bar P)$ (from now on we
abbreviate $\R_w (\pi_1 X;SO(3))$ to $\R_w (X;SO(3))$ etc). The
following result lifts this correspondence to the level of projectively flat
connections. It is proved in \cite{ruberman-saveliev:casson}.

\begin{proposition}\label{P:hol}
Let $P$ be an admissible $U(2)$--bundle over a homology $4$--torus $X$, and
$c$ a $2$--cocycle representing $w = w_2 (\bar P)$. Then the holonomy
correspondence defines a map $\phi: \M(P)\to \P_c(X; SU(2))$ which is
an $H^1(X;\Z_2)$--equivariant bijection.
\end{proposition}

In particular, we have the following commutative diagram whose horizontal
arrows are bijections
\[
\begin{CD}
\M(P)       @> \phi     >>  \P_c (X;SU(2)) \\
@VV\pi V                     @VV\pi V \\
\M(\bar P)  @> \bar\phi >>  \R_w (X;SO(3))
\end{CD}
\]
Here, $\pi: \P_c (X;SU(2)) \to \R_w (X;SO(3))$ is the natural projection
map, see Section \ref{S:algebra}.

An argument similar to that for representation varieties shows that Zariski
tangent space to $\P_c (X;SU(2))$ at a projective representation $\rho:
\pi_1(X)\to SU(2)$ equals $H^1 (X;\ad\rho)$ where $\ad\rho: \pi_1(X)\to SU(2)
\to SO(3)$ is a representation. It is identified as usual with the tangent
space to $\M(P)$ at the corresponding projectively flat connection.

A point $\rho\in\P_c (X;SU(2))$ is called \emph{non--degenerate} if $H^1(X;
\ad\rho) = 0$; the space $\P_c (X; SU(2))$ is called \emph{non--degenerate}
if all of its points are non--degenerate. Since $p_1(\bar P) = 0$ the formal
dimension of $\M(P)$ equals zero, and so $H^1(X;\ad\rho)$ and $H^2_+
(X;\ad\rho)$ have the same dimension. In particular, $\M(P)$ is
non--degenerate if and only if $\P_c (X;SU(2))$ is, compare with Section
\ref{S:proj-asd}. Because of the identification of Proposition \ref{P:hol},
the invariant $\bar\lambda(X,P)$ in this situation can be computed as
signed count of points in $\P_c (X;SU(2))$ with $[c] = w_2(\bar P)$.


\section{The four--orbits}\label{S:four-orbits}
Let $X$ be a homology $4$--torus. According to the action of $H^1(X;\Z_2) =
(\Z_2)^4$ the space $\P_c(X;SU(2))$ splits into orbits of possible orders
one, two, four, eight, and sixteen. In this section we study the
\emph{four--orbits} (orbits with four elements, or those with stabilizer
$\ZZ$).  The results of this section are quite parallel to those in Section
4 of~\cite{ruberman-saveliev:casson}, and the proofs are only sketched.


\subsection{The four--orbits and invariant $\bar\lambda$}

Consider a subgroup of $SO(3)$ that is isomorphic to $\ZZ$. Such a subgroup
is generated by $\pi$--rotations about two perpendicular axes in
$\mathbb R^3$, and any two such subgroups are conjugate to each other in
$SO(3)$. Hence the following definition makes sense. Given $w\in H^2
(\pi_1 X; \Z_2)$, define $\R_w (X;\ZZ)$ to be the (finite) subspace of
$\R_w (X;SO(3))$ consisting of the $SO(3)$ conjugacy classes of
representations $\alpha: \pi_1(X)\to SO(3)$ which factor through $\ZZ
\subset SO(3)$, compare with Section \ref{S:algebra}.

\begin{proposition}
Let $[c] = w$ be a non-trivial class in $H^2(\pi_1 X;\Z_2)$. Then the map
$\pi: \P_c(X;SU(2))\to \R_w(X;SO(3))$ establishes a bijective correspondence
between the set of four--orbits in $\P_c(X;SU(2))$ and the set $\R_w(X;\ZZ)$.
\end{proposition}

\begin{proof}
Suppose that the conjugacy class of a projective representation $\rho:
\pi_1(X)\to SU(2)$ is fixed by a subgroup $\ZZ$ of $H^1(X;\Z_2)$ generated
by homomorphisms $\alpha$, $\beta: \pi_1(X) \to \Z_2$.  As
in~\cite[Proposition 4.1]{ruberman-saveliev:casson}, $\rho$ factors
through a copy of the quaternion $8$-group $Q = \{\,\pm 1, \pm i,\pm
j, \pm k\,\}$ and therefore its associated $SO(3)$ representation $\ad\rho$
factors through a copy of $\ZZ\subset SO(3)$.

To complete the proof, we only need to show that the orbit of $\rho$
consists of exactly four points. Let $\gamma$, $\delta$ be vectors in
$H^1(X;\Z_2)$ completing $\alpha$, $\beta$ to a basis. Then $\rho$,
$\rho^{\gamma}$, $\rho^{\delta}$ and $\rho^{\gamma\delta}$ lie in the
same $H^1(X;\Z_2)$--orbit but are not conjugate. The latter can be seen as
follows: suppose there exists an element $w\in SU(2)$ such that
$\gamma(x)\rho(x)
= w\,\rho\,w^{-1}$ then $w = \pm k$ and $\alpha(x)\beta(x)\gamma(x)\rho(x)
= (ijk)\,\rho(x)\, (ijk)^{-1} = \rho(x)$ for all $x$, a contradiction.
\end{proof}

\begin{remark}\label{R:one-orbits}
The above proof also shows that no point of $\P_c (X;SU(2))$ with
$[c]\ne 0$ is fixed by more than two basis elements in $H^1(X;\Z_2)$ so
that $\P_c(X;\allowbreak SU(2))$ has no orbits of orders one or two.
\end{remark}


\subsection{The number of four--orbits}

Our next goal is to find a formula for the number of points in $\R_w(X;\ZZ)$
modulo 2.

\begin{proposition}\label{P:cup}
If $X$ is an odd homology $4$--torus then for any non-zero $w\in H^2(\pi_1 X;
\Z_2)$ such that $w\cup w \equiv 0\mod 4$ we have $\#\,\R_w (X;\ZZ) \equiv
1\pmod 2$.
\end{proposition}

\begin{proof}
Recall that any two subgroups of $SO(3)$ that are isomorphic to $\ZZ$ are
conjugate, and that moreover any automorphism of such a subgroup is
realized by conjugation by an element of $SO(3)$. Let us fix a subgroup
$\ZZ$ and a basis in it.

Since $\ZZ$ is abelian, every $\alpha\in\R_w (X;\ZZ)$ factors through a
homomorphism $H_1(X;\Z)\to \ZZ$. The two components of this homomorphism
determine elements $\beta,\gamma \in\Hom(H_1(X;\Z);\Z_2)\cong H^1(X;\Z_2)$.
It is straightforward to see that the $SO(3)$ representation $\alpha$ may
be recovered from $\beta$ and $\gamma$ via the formula $\alpha\cong \beta
\oplus \gamma \oplus \det(\beta \oplus \gamma)$. This establishes a
bijective correspondence
\begin{equation}\label{E:lambda}
\R(X;\ZZ) \to \Lambda^2_0 H^1(X;\Z_2),
\end{equation}
where $\R(X;\ZZ)$ is union of $\R_w (X;\ZZ)$ over all possible $w$ and
$\Lambda^2_0 H^1(X;\Z_2)$ is the subset of decomposable elements of
$\Lambda^2 H^1(X;\Z_2)$. Since $H_1(X;\Z)$ is torsion free, any element
in $H^1(X;\Z_2)$ is the mod 2 reduction of a class in $H^1(X;\Z)$. It
follows that the cup product of any element $a\in H^1(X;\Z_2)$ with
itself is $0$. We compute
\begin{equation}\label{whitney}
\begin{split}
w_2(\alpha) & = w_1(\beta) \cup w_1(\gamma) +  w_1(\beta) \cup
w_1(\det(\beta \oplus \gamma)) + w_1(\gamma) \cup w_1(\det(\beta
\oplus \gamma)) \\
& = w_1(\beta) \cup w_1(\gamma) + w_1(\beta) \cup (w_1(\beta) + w_1(\gamma)) +
w_1(\gamma) \cup (w_1(\beta) + w_1(\gamma))\\
&= w_1(\beta) \cup w_1(\gamma).
\end{split}
\end{equation}
Since $w_1(\beta) = \beta$ and $w_1(\gamma) = \gamma$, this shows that
$w_2(\alpha)$ belongs to the image of $\Lambda_0^2 H^1 (X;\Z_2)$ under the
map (\ref{E:cup}).

The result now follows by combining (\ref{E:cup}) and (\ref{E:lambda})\,:
Corollary \ref{C:cup} implies that the space $\R_w (X;\ZZ)$ consists of
exactly one point for every choice of $0\ne w\in H^2(\pi_1 X;\Z_2)$ such
that $w\cup w = 0\pmod 4$.
\end{proof}

\begin{proposition}
If $X$ is an even homology $4$--torus then for any non-zero $w\in H^2
(\pi_1 X; \Z_2)$ such that $w\cup w\equiv 0\pmod 4$ and $w\cup\xi\ne 0$
for some $\xi\in H^1(X;\Z_2)$ the space $\R_w (X;\ZZ)$ is empty.
\end{proposition}

\begin{proof}
This follows by combining (\ref{E:cup}) and (\ref{E:lambda}) and using the
following observation. Suppose that $w\in H^2 (X; \Z_2)$ is of the form $w
= \beta \cup \gamma$ then  $\omega \cup \xi = \beta \cup \gamma \cup \xi =
0$ for all $\xi \in H^1 (X;\Z_2)$\,: if not, then by Poincar\'e duality,
there exists $\eta\in H^1 (X;\Z_2)$ such that $\beta \cup \gamma \cup \xi
\cup \eta\ne 0\pmod 2$, which contradicts the fact that $X$ is an even
homology torus.
\end{proof}

\begin{corollary}\label{C:det}
Let $P$ be an admissible $U(2)$--bundle over a homology $4$--torus $X$ and let
$w = w_2 (\bar P)$ be in the image of (\ref{E:iota}). Then $\#\,\R_w (X;\ZZ)
\equiv \det X \pmod 2$.
\end{corollary}


\subsection{Non--degeneracy of the four--orbits}\label{S:4-nondeg}
We wish to use Corollary \ref{C:det} to calculate the contribution of the
four--orbits into $\bar\lambda(X,P)\pmod 2$.  In order to do that, we need
to check the non--degeneracy condition for such orbits in the case when
$X$ is an odd homology $4$--torus. This is similar to the non--degeneracy
statement about two--orbits in~\cite{ruberman-saveliev:casson}.

\begin{proposition}\label{P:4-nondeg}
If $X$ is an odd homology $4$--torus then $H^1(X;\ad\rho) = 0$ for any
projective representation $\rho: \pi_1(X)\to SU(2)$ such that
$\ad\rho \in\R_w (X;\ZZ)$ with $w\ne 0$ and $w\cup w = 0\pmod 4$.
\end{proposition}

As in the proof of~\cite[Proposition 4.4]{ruberman-saveliev:casson},
the cohomology group in question is related to the usual cohomology
of $X$ via a Gysin sequence.  The non-triviality of the cup product
determines enough about the maps to lead to the non-vanishing result.
We refer the reader to~\cite{ruberman-saveliev:casson} for all the
details.


\section{The eight--orbits}\label{S:eight}
Let $X$ be a $\Z[\Z]$--homology $4$--torus and $P$ an admissible $U(2)$--bundle
over $X$. The moduli space $\M(P)$ need not be non--degenerate at the
eight--orbits as it was at the four--orbits and therefore has to be perturbed.
Doing so equivariantly requires a detailed knowledge of the geometry of
$\M(P)$ at the eight--orbits, which is the topic of the first part of this
section. We proceed by perturbing $\M(P)$ at the eight--orbits, first in
the tangential and then in the normal direction.


\subsection{Connections with stabilizer $\Z_2$}
Let us fix $0 \ne \chi\in H^1(X;\Z_2)$ and consider the action $\chi:
\B^*(P)\to \B^*(P)$ induced on the irreducible part of $\B(P)$ by
(\ref{E:chi}). The fixed points of this action will be denoted by
$\B^{\chi}(P)$, and we begin by describing $\B^{\chi}(P)$ as a
moduli space in its own right.  To this end we make the following definitions.

Let $u: P \to P$ be a $U(2)$ gauge transformation realizing $\chi\in
H^1 (X; \Z_2)$ and let $h = \det u$. Then we denote by $\A^u(P)$ the
subset of $\A(P)$ consisting of connections $A$ such that $u^* A -
1/2\,h^{-1} dh\,\Id = A$, and by $\G^u (P)$ the subgroup of $\G(P)$
consisting of gauge transformations $g \in \G(P)$ having the property
that $g u = \pm\,u g$.

\begin{proposition}\label{P:chi}
Suppose that $\B^{\chi}(P)$ is non--empty. If $u: P \to P$ realizes $\chi
\in H^1 (X;\Z_2)$ then $\B^{\chi}(P) = \A^u(P)/\,\G^u(P)$.  Moreover,
$\bar u^2 = \Id$.
\end{proposition}

\begin{proof}
Define a map $\A^u(P)\to \B^{\chi} (P)$ by $A \to [A]$. This map is
surjective: given $[A]\in \B^{\chi}(P)$, there exists a lift $g$ of $\chi$
and a gauge transformation $k \in \G(P)$ such that $g^* A - 1/2\,\det g
\cdot d\,(\det g)^{-1}\Id = k^* A$. Let $v = k^{-1} g$ then $v^* A - 1/2\,
\det v\cdot d\,(\det v)^{-1}\Id = A$ since $\det g = \det v$ and $k$
commutes with $d (\det g^{-1})\Id$. Since both $u$ and $v$ lift
$\chi$, there exists an automorphism $\gamma = c\,\Id$
with $c \in U(1)$ such that $\gamma v = w^{-1} u\,w$ for some $w\in \G(P)$.
A calculation similar to (\ref{E:calc}) shows that
\begin{multline}\notag
v^* A - 1/2\,\det v\cdot d\,(\det v)^{-1}\,\Id
= (\gamma v)^* A - 1/2\,\det \gamma v\cdot d\,(\det \gamma v)^{-1}\Id \\
= (w^{-1} u\,w)^* A - 1/2\,\det u\cdot d\,(\det u)^{-1}\Id = A.
\end{multline}
Therefore, $u^* w^* A - 1/2\,h\,dh^{-1}\Id = w^* A$, which implies that
$[A] = [w^* A]$ with $w^* A \in \A^u(P)$.

The map $\A^u(P)\to \B^{\chi}(P)$ becomes injective after we factor
$\A^u(P)$ out by the gauge group $\G^u(P)$. This can be seen as follows.
Suppose that $A, B \in \A^u(P)$ are such that $[A] = [B]$,  then there is
$g \in \G(P)$ such that $B = g^* A$ and $u^* g^* A - 1/2\,h\,dh^{-1}\Id =
g^* A$.  Therefore,  in addition to $u^* A - 1/2\,h\,dh^{-1}\Id = A$,  we
also have $(g^{-1} u g)^* A - 1/2\,h\,dh^{-1}\Id = A$.  This implies that
$(g^{-1} u g)^* A = u^* A$ and that $(u^{-1} g^{-1} u g)^* A = A$. Since
$\det (u^{-1} g^{-1} u g) = 1$ and $A$ is irreducible,  we conclude that
$u^{-1} g^{-1} u g = \pm 1$.

The equality $(u^2)^*A - h\,dh^{-1}\Id = A$ obtained by applying $u^*$ to
$A$ twice is equivalent to $(h^{-1} u^2)^* A = A$ which is seen from the
following calculation\,:  $(h^{-1} u^2)^* A = h^{-1}dh\,\Id + (u^2)^* A =
h^{-1}dh\,\Id + h\,dh^{-1}\Id + A = A$. Since $\det (h^{-1} u^2) = h^{-2}
h^2 = 1$ and $A$ is irreducible, we conclude that $h^{-1} u^2 = \pm 1$ and
therefore $\bar u^2 = 1$.
\end{proof}


\subsection{The local structure of $\M(P)$ near eight--orbits}
It follows from Corollary \ref{C:chi} that the action $\chi: \B(P)\to\B(P)$
restricts to a well--defined action $\chi: \M(P) \to \M(P)$,  whose fixed
point set will be denoted by $\M^{\chi}(P) \subset \B^{\chi}(P)$.  We wish
to investigate the non--degeneracy condition for $\M(P)$ at points
$[A]\in\M^{\chi}(P)$.
If $\M(P)$ fails to be non--degenerate at $[A] \in \M^{\chi}(P)$,  this may
be due to the nonvanishing of the obstruction (lying in the second
cohomology of the deformation complex) to infinitesimal deformations in the
direction of $\M^{\chi}(P)$ (tangential deformations) and/or in the
direction normal to $\M^{\chi}(P)\subset \M(P)$.

More precisely, according to Proposition \ref{P:chi}, one can choose a
gauge transformation $u: P \to P$ realizing $\chi$ so that every point in
$\M^{\chi}(P)$ is the $\G^u(P)$ equivalence class of a connection $A$
having the property $u^*A - 1/2\,h\,dh^{-1}\Id = A$ with $h = \det u$.
The isomorphism (\ref{D:iso}) of deformation complexes at $[A]$ and
$\chi[A] = [A]$ then becomes an order two automorphism $(\ad\bar u)^*$ of
the deformation complex at $[A]$,  thus splitting it into a direct sum of
two elliptic complexes according to the $(\pm 1)$--eigenvalues of the
operator $(\ad\bar u)^*$. The complex
\begin{equation}\label{D:plus}
\begin{CD}
\Omega^0 (X;\ad\bar P)^+ @> d_{\bar A} >> \Omega^1 (X;\ad\bar P)^+
@> d^+_{\bar A} >> \Omega^2_+ (X;\ad\bar P)^+
\end{CD}
\end{equation}
describes tangential deformations, while normal deformations are
described by the complex
\begin{equation}\label{D:minus}
\begin{CD}
\Omega^0 (X;\ad\bar P)^- @> d_{\bar A} >> \Omega^1 (X;\ad\bar P)^-
@> d^+_{\bar A} >> \Omega^2_+ (X;\ad\bar P)^-.
\end{CD}
\end{equation}
The cohomology of the above two complexes will be denoted by
$H^0(X;\ad A)^{\pm}$, $H^1(X;\ad A)^{\pm}$, and $H^2_+(X;\ad A)^{\pm}$,
respectively.

\begin{lemma}
Both complexes (\ref{D:plus}) and (\ref{D:minus}) have zero index.
\end{lemma}

\begin{proof}
Since the index of the deformation complex (\ref{D:full}) is zero,   it is
sufficient to show that the complex (\ref{D:plus}) has zero index. Observe
that the $(+1)$--eigenspace $\ad\bar P^{\,+}$  of the operator $\ad\bar u:
\ad\bar P \to \ad\bar P$ is a $\Z_2$--bundle with fiber $\mathbb R$ and
$w_1 (\ad\bar P^+) = \chi$, and that $\Omega^k (X;\ad\bar P)^+ = \Omega^k
(X;\ad\bar P^{\,+})$ for all $k$. Therefore, the complex (\ref{D:plus}) is
isomorphic to the deformation complex of the moduli space of ASD
connections on $\ad\bar P^{\,+}$, and the latter has index
$-(1 - b_1 + b_2^+)\,(X) = 0$ by the index theorem.
\end{proof}


\subsection{Tangential perturbations}\label{S:tangent}
Let $\sigma: \B(P) \to \E(P)$ be an equivariant admissible perturbation,
and $\M_{\sigma} (P)$ the corresponding moduli space of perturbed
projectively ASD connections. Every $0\ne \chi \in H^1(X;\Z_2)$ defines
a map $\chi: \M_{\sigma} (P) \to \M_{\sigma} (P)$ whose fixed point set
will be denoted by $\M^{\chi}_{\sigma} (P)$. Let us choose a gauge
transformation $u: P \to P$ as in Proposition \ref{P:chi} then the
deformation complex of $\M^{\chi}_{\sigma}(P)$ at $[A]$,
\begin{equation}\label{D:pert+}
\begin{CD}
\Omega^0 (X;\ad\bar P)^+ @> d_{\bar A} >> \Omega^1 (X;\ad\bar P)^+
@> d^+_{\bar A} + D_A \sigma >> \Omega^2_+ (X;\ad\bar P)^+,
\end{CD}
\end{equation}
is the invariant part of the complex (\ref{D:pert}) with respect to the
action of $(\ad\bar u)^*$. We say that $\M^{\chi}_{\sigma}(P)$ is
\emph{non--degenerate} at $[A]$ if the second cohomology of the complex
(\ref{D:pert+}) vanishes (observe that, in general, this condition may be
weaker than the non--degeneracy condition for the full moduli space
$\M_{\sigma}(P)$ at $[A]$). If $\sigma = 0$, this condition is equivalent
to the vanishing of the second cohomology group $H^2_+ (X;\ad A)^+$ of the
complex (\ref{D:plus}). We call $\M^{\chi}_{\sigma}(P)$
\emph{non--degenerate} if it is non--degenerate at all $[A] \in
\M^{\chi}_{\sigma}(P)$. Our goal in this section will be to find enough
equivariant admissible perturbations $\sigma$ to ensure that
$\M^{\chi}_{\sigma}(P)$ is non--degenerate for a generic $\sigma$.

We begin by restating the above perturbation problem in slightly
different terms. Proposition \ref{P:equiv} implies that, for every
$0 \ne \chi \in H^1(X;\Z_2)$, we have a commutative diagram
\[
\begin{CD}
\B(P) @> \sigma >> \E(P) \\
@VV \chi V    @VV \chi V \\
\B(P) @> \sigma >> \E(P).
\end{CD}
\]
In particular, the restriction of $\sigma$  onto $\B^{\chi}(P)$ defines a
map $\sigma^{\chi}: \B^{\chi}(P) \to \E^{\chi}(P)$,  where $\E^{\chi}(P)$
stands for the fixed point set of $\chi: \E(P)\to \E(P)$. For a choice of
$u: P \to P$ as above, $\E^{\chi}(P)$ is isomorphic to the bundle
\[
\E^{\chi}(P) = \A^u(P)\, \times_{\,\G^u(P)}\, \Omega^2_+ (X;\ad\bar P)^+
\]
over $\B^{\chi}(P)$, and $\sigma^{\chi}$ is a section of this bundle. The
restriction of the section $F_+: \B(P) \to \E(P)$ onto $\B^{\chi}(P)$
defines a section $F_+^{\chi}: \B^{\chi}(P) \to \E^{\chi}(P)$ in a similar
fashion. The zero set of $\zeta^{\chi}_{\sigma} = F_+^{\chi} +
\sigma^{\chi}: \B^{\chi}(P) \to \E^{\chi}(P)$ is precisely the moduli
space $\M^{\chi}_{\sigma}(P)$, only now viewed as obtained by perturbing
$\M^{\chi}(P)$ using perturbation $\sigma^{\chi}$.

This reformulation makes the tangential perturbation problem similar to the
non--equivariant perturbation problem studied in Section \ref{S:adm-pert}.
The key reduction step is provided by the following lemma.

\begin{lemma}\label{L:lifts}
Let $p:\tilde X\to X$ be the regular 16--fold covering of $X$ corresponding
to the mod 2 abelianization homomorphism  $\pi_1 X \to H_1 (X;\Z_2)$.   For
every projectively flat connection  $A$ on $P \to X$,  denote by $\tilde A$
the pull back connection over $\tilde X$. Let $\Stab(A)$ denote the
stabilizer of $[A]$ in $H^1(X;\Z_2)$ then
\begin{enumerate}
\item[\rm(a)] $\Stab(A) = 1$ if and only if $\tilde A$ is irreducible,
\item[\rm(b)] $\Stab(A) = \Z_2$ if and only if $\tilde A$ is reducible
non--central, and
\item[\rm(c)] $\Stab(A) = \Z_2\,\oplus\,\Z_2$ if and only if $\tilde A$
is central.
\end{enumerate}
No other stabilizers $\Stab(A)$ may occur.
\end{lemma}

\begin{proof}
A version of this statement for the homology $3$--tori is proved in
\cite[Lemma 5.6]{ruberman-saveliev:casson}; the proof works exactly the
same in the $4$--dimensional case.
\end{proof}

\begin{lemma}\label{L:ab8}
Let $A$ be a projectively flat connection with stabilizer $\Z_2$ in
$H^1(X;\Z_2)$. Then there exists a collection of disjoint mod-2 trivial
loops $\psi_k$ and a collection of self--dual $2$--forms $\nu_k$, each
supported in its respective $\psi_k (S^1\times D^3)$, such that the forms
$\Pi\hol_A(\psi_k)\,\otimes\,\nu_k$ generate the vector space
$H^2_+\,(X;\ad A)^+$.
\end{lemma}

\begin{proof}
Fix a harmonic lift $H^2_+(X;\ad A)^+ \subset \Omega^2_+(X;\ad\bar P)^+$
of $H^2_+(X;\ad A)^+$.  It consists of $2$--forms $\omega$ with coefficients
in $\ad\bar P$ such that $\omega$ is self--dual, $d_{\bar A} \omega = 0$,
and $(\ad\bar u)^*\omega = \omega$. Observe that, for a mod-2 trivial loop
$\psi_k$, the equality $u^*A - 1/2\,h\,dh^{-1}\Id = A$ with $h = \det u$
implies that
\begin{multline}\notag
\Pi\hol_A(\psi_k)\,\otimes\,\nu_k =
\Pi\hol_{u^* A - 1/2\,h\,dh^{-1}\Id}(\psi_k)\,\otimes\,\nu_k \\
= (\ad\bar u)^*\; (\,\Pi\hol_A(\psi_k)\,\otimes\,\nu_k\,),
\end{multline}
as in the proof of Proposition \ref{P:equiv}, hence the forms
$\Pi\hol_A (\psi_k)\,\otimes\,\nu_k$ belong to $\Omega^2_+(X;\ad\bar P)^+$.
We need to construct sufficiently many of these forms so that the subspace
of $\Omega^2_+ (X;\ad\bar P)^+$ that they span projects orthogonally
\emph{onto} $H^2_+ (X;\ad A)^+$.

There exist finitely many points $x_1,\ldots, x_m$ in $X$ such that the
evaluation map
\[
\ev: H^2_+(X;\ad A)^+ \longrightarrow\; \bigoplus_{p = 1}^m\;\Lambda^2_+
(\ad\bar P)^+_{x_p}
\]
given by $\ev(\omega) = (\omega(x_1),\ldots,\omega(x_m))$ is injective. Here,
$(\ad\bar P)^+_{x_p}$ is the $(+1)$--eigenspace of $(\ad\bar u)_{x_p}$ acting
on $(\ad\bar P)_{x_p}$. Since the latter eigenspace is one--dimen\-sional,
each of the spaces $\Lambda^2_+(\ad\bar P)^+_{x_p}$ has real dimension three.

Let $\rho: \pi_1 X \to SU(2)$ be the holonomy representation of $A$. For each
of the points $x_p$, consider the possible holonomies $\rho(\psi)$ around
mod-2 trivial loops $\psi$ based at $x_p$. According to Lemma \ref{L:lifts},
the restriction of $\rho$ to the subgroup $\pi_1 \tilde X \subset \pi_1 X$
generated by mod-2 trivial loops is non--central, so there exists a mod-2
trivial loop $\psi_p$ such that $\Pi\hol_A (\psi_p) \in (\ad\bar P)^+_{x_p}$
is not zero.

Let $\omega_{i,p}$, $i = 1, 2, 3$, be a basis for $(\Lambda^2_+)_{x_p}$ then
the tensor products $\Pi\hol_A (\psi_p)\,\otimes\,\omega_{i,p}$ span the
entire space
\[
\bigoplus_{p = 1}^m\;\Lambda^2_+(\ad\bar P)^+_{x_p}.
\]
For each index $p$, choose a small ball around $x_p$ over which the forms in
$\Omega^2_+(X;\ad\bar P)^+$ have small variations, and three distinct points
inside it. Label the three points near $x_p$ by $x_{i,p}$. Then we can
choose mod-2 trivial loops $\psi_{i,p}$ based at $x_{i,p}$ whose holonomy
is close to $\hol_A (\psi_p)$, and also $2$--forms $\nu_{i,p}$ supported in
small balls around $x_{i,p}$ and close to the respective multiples of
$\omega_{i,p}$ in a local trivialization. The resulting sections $\Pi\hol_A
(\psi_{\,i,p})\,\otimes\,\nu_{i,p}$ have disjoint supports and, if the
approximations in the above construction are made sufficiently fine, no
non-zero element of $H^2_+(X;\ad A)^+$ can be orthogonal to all of them.
\end{proof}

Observe that, for every $0 \ne \chi \in H^1(X; \Z_2)$, the moduli space
$\M^{\chi}(P)$ is compact as it is the fixed point set of an involution
acting on the compact space $\M(P)$. In particular, we can suppose that
the mod-2 trivial loops $\psi_k$ and forms $\nu_k$ can be chosen so that
the conclusion of Lemma \ref{L:ab8} holds for all $[A] \in \M^{\chi}(P)$
simultaneously. We fix such a choice and define $\sigma$ by the formula
(\ref{E:sigma'}).

\begin{proposition}
There exists a real number $r > 0$ such that, with the choice of $\psi_k$
and $\nu_k$ as above, the moduli space $\M^{\chi}_{\sigma}(P)$ is
non--degenerate for a generic $\bar\ep = (\ep_1,\ldots,\ep_N)$ inside the
ball $\|\bar\ep\| < r$ in $\mathbb R^N$.
\end{proposition}

\begin{proof}
Let us pull the bundle $\E^{\chi}(P)$ back to a bundle $\pi^*\,\E^{\chi}(P)
\to \B^{\chi}(P) \times \mathbb R^N$ and consider the ``universal section''
$\Psi^{\chi}$ obtained by restricting the section (\ref{E:psi}) to
$\B^{\chi}(P) \times \mathbb R^N$. Since the forms $\Pi\hol_A(\psi_k)\,
\otimes\,\nu_k$ generate $H^2_+(X;\ad A)^+$, we conclude that $\Psi^{\chi}$
is transversal at $\M^{\chi}(P)\times\{0\}$ to the zero section of the
above bundle. The rest of the argument is identical to the proof of
Proposition \ref{P:gen}.
\end{proof}

One can use the holonomy correspondence of  Section \ref{S:hol} to conclude
that $\M^{\chi}(P)$ for each $0 \ne \chi \in H^1(X; \Z_2)$ contains at most
one four--orbit. Since the four--orbits are non--degenerate, they are
isolated in each of the $\M^{\chi}(P)$ and hence the entire moduli space of
eight--orbits is compact. The above argument then shows that, after a
generic equivariant admissible perturbation, the moduli space of
eight--orbits is non--degenerate.


\subsection{Normal perturbations}\label{S:tau}
The non--degeneracy of the moduli space of eight--orbits achieved in the
previous section by the means of a perturbation $\sigma$ implies that the
moduli space $\M_{\sigma}(P)$ has finitely many eight--orbits, and it is
non--degenerate at each of them in the tangential direction. However,
$\M_{\sigma}(P)$ may fail to be non--degenerate in the normal direction,
in which case another perturbation is needed. Such a perturbation, if it
is small enough, will not spoil the tangential non--degeneracy of
$\M_{\sigma}(P)$. This observation, together with the fact that there are
only finitely many eight--orbits in $\M_{\sigma}(P)$, allows us to perturb
at each of them separately. Furthermore, to save some notation, we will
assume that the \emph{unperturbed} moduli space $\M^{\chi}(P)$ is
non--degenerate; the argument below remains valid after a small
equivariant perturbation.

First we wish to extend the class of admissible perturbations, following
ideas of Donaldson \cite{donaldson} and Furuta 
\cite{furuta:perturbation}. Assume that
$\M^{\chi}(P)$ is non--degenerate at $[A]$ and choose an open neighborhood
$\U([A]) \subset \B(P)$ which does not contain any other points of $\M(P)$
with non--trivial stabilizer. The slice theorem then implies that, if
$\U([A])$ is small enough, one can find a real number $r > 0$ such that
every $B\in \A(P)$ with $[B]\in\U([A])$ can be written as $B = g^*(A + a)$
for a unique (up to sign) $g\in\G(P)$ and a unique $a\in\ker d^*_{\bar A}$
such that $\|a\| = \|a\|_{L^2_l} < r$.

Choose a bounded linear operator  $\phi: \ker d^*_{\bar A} \to \Omega^2_+
(X; \ad\bar P)$, where both spaces are completed in Sobolev $L^2_l$ norms,
and let $\beta: \mathbb R_+ \to [0,1]$  be a smooth cut off function such
that $\beta(t) = 1$ for $t \le r/2$ and $\beta(t) = 0$ for $t \ge r$. The
formula
\begin{equation}\label{E:tau-def}
\tau_{[A]}([B]) = [\,A + a,\;\beta(\|a\|)\,\phi(a)\,],\quad [B]\in \U([A]),
\end{equation}
extended to be zero on the complement of $\U([A])$, defines a section
$\tau_{[A]}: \B(P) \to \E(P)$.

Next we address equivariance properties of the above section.  Remember
that $[A] \in \M^{\chi}(P)$  for some $0 \ne \chi \in H^1 (X;\Z_2)$, so
we first handle the equivariance of $\tau_{[A]}:  \B(P) \to \E(P)$ with
respect to the action of $\chi$. Fix a representative $A$ and use
Proposition \ref{P:chi} to realize $\chi$ by a gauge transformation $u:
P \to P$ such that $u^* A - 1/2\,h\,dh^{-1}\,\Id = A$.

\begin{lemma}
The perturbation $\tau_{[A]}$ is equivariant with respect to the action
of $\chi$ if and only if the following diagram commutes
\begin{equation}\label{D:phi}
\begin{CD}
\ker d^*_{\bar A} @> \phi >> \Omega^2_+ (X; \ad\bar P) \\
@VV (\ad\bar u)^* V  @VV (\ad\bar u)^* V \\
\ker d^*_{\bar A} @> \phi >> \Omega^2_+ (X; \ad\bar P).
\end{CD}
\end{equation}
\end{lemma}

\begin{proof}
Observe that the operator $(\ad\bar u)^*$  is an isometry on the Sobolev
$L^2_l$ completions of $\ker d^*_{\bar A}$ and $\Omega^2_+(X;\ad\bar P)$.
The statement now follows by a straightforward calculation using the
formula (\ref{E:tau-def}).
\end{proof}

Thus, for any $[A] \in \M^{\chi}(P)$, the equivariance of $\tau_{[A]}$
with respect to the action of $\chi$ can be achieved by making a proper
choice of $\phi$ in its definition, see Lemma \ref{L:tau_k} below. The
residual action of $H^1 (X; \Z_2)/\chi = (\Z_2)^3$ can then be used to
spread the perturbation $\tau_{[A]}$ to a small neighborhood of the
orbit of $[A]$. More explicitly, if $\eta \in H^1(X;\Z_2)/\chi$, then
$\U(\eta[A]) = \eta\,\U([A])$ is a neighborhood of $\eta[A]$, and we
can assume that all of the sets $\U(\eta[A])$ have disjoint
closures. Then $\tau_{\eta [A]} = \eta\,\tau_{[A]}\,\eta^{-1}$ is a
perturbation that is equivariant with respect to the stabilizer
$\chi\eta^{-1}$ of $\eta[A]$. Hence
\begin{equation}\label{E:tau-def'}
\tau = \sum\; \tau_{\eta [A]}: \B(P) \to \E(P)
\end{equation}
is a section which is equivariant with respect to the action of the full
group  $H^1 (X; \Z_2)$ and which is supported in a small neighborhood of
the eight--orbit through $[A]$.

We define new admissible perturbations $\tau: \B(P)\to\E(P)$ as finite
linear combinations of perturbations defined in (\ref{E:tau-def'}).
Each section $\tau$ is smooth and compact since it factors through the
embedding $\E'(P) \to \E(P)$, where $\E'(P)$ is the bundle (\ref{E:bndl})
whose fiber is completed in Sobolev $L^2_l$ norm, and since the inclusion
$L^2_l \to L^2_{l-1}$ is compact, compare with Section \ref{S:adm-pert}.
Therefore, the section $\zeta_{\tau} = F_+ + \tau: \B(P) \to \E(P)$ is
Fredholm, and the perturbed moduli space $\M_{\tau}(P) =
(\zeta_{\tau})^{-1}(0)$ has formal dimension zero.

\begin{proposition}
For any choice of perturbation $\tau$ as above, the moduli space
$\M_{\tau}(P)$ is compact.
\end{proposition}

\begin{proof}
The proof of Proposition \ref{P:compact} goes through as before to give
the above statement because the curvatures $F_+ (\bar A)$ are uniformly
bounded in $L^{\infty}$ for all $[A] \in \M_{\tau}(P)$, see Furuta
\cite[Proposition 3.1]{furuta:perturbation}.
\end{proof}

Our next step is to show that the operators $\phi$ we use to define $\tau$
can be chosen so that $\M_{\tau} (P)$ is acted upon by $H^1 (X; \Z_2)$ and
is non--degenerate at the eight--orbits of this action.

\begin{lemma}\label{L:tau_k}
For any $[A] \in \M^{\chi} (P)$, there exist finitely many bounded linear
operators $\phi_k: \ker d^*_{\bar A} \to \Omega^2_+ (X; \ad\bar P)$ which
make the diagram (\ref{D:phi}) commute and have the following properties.
Let $\tau_k$ be the perturbation (\ref{E:tau-def'}) defined using $\phi_k$
and $D_A\tau_k: H^1 (X;\ad A)^- \to H^2_+ (X;\ad A)^-$ the linear map
obtained by restricting the derivative of $\tau_k$ on $H^1 (X; \ad A)^-
\subset \Omega^1 (X, \ad\bar P)^-$ and projecting onto $H^2_+ (X;\ad A)^-
\subset \Omega^2_+ (X; \ad\bar P)^-$.   Then
the maps $D_A \tau_k$ span $\Hom\,(H^1(X;\ad A)^-, H^2_+ (X;\ad A)^-)$.
\end{lemma}

\begin{proof}
The derivative of $\tau_k$ at $A$ is the operator $\phi_k$. Let us choose
a basis $\psi_k$ for $\Hom\,(H^1(X;\ad A)^-,H^2_+(X;\ad A)^-)$ and extend
$\psi_k$ by composing with the projection $\ker d_{\bar A}^* \to H^1 (X;
\ad A)^-$ to obtain a family of bounded linear operators $\phi_k: \ker
d^*_{\bar A} \to \Omega^2_+ (X, \ad\bar P)$. The action of $(\ad\bar u)^*$
on $H^1 (X; \ad A)^-$ and $H^2_+ (X;\ad A)^-$ is by minus identity,  hence
commutativity of (\ref{D:phi}) is immediate from the linearity of $\phi_k$.
\end{proof}

\begin{remark}\label{R:pert}
If we start with an already perturbed moduli space $\M_{\sigma}(P)$, where
$\sigma$ is any linear combination of equivariant perturbations of the
types (\ref{E:sigma'}) and (\ref{E:tau-def'}), the above lemma can easily
be amended to provide perturbations $\tau_k$ whose derivatives span $\Hom\,
(H^1_{\sigma} (X;\ad A)^-, H^2_{+,\sigma} (X;\ad A)^-)$, where
$H^1_{\sigma}(X;\ad A)^-$ and $H^2_{+, \sigma} (X;\ad A)^-$ are the
cohomology groups of the complex
\[
\begin{CD}
\Omega^0 (X;\ad\bar P)^- @> d_{\bar A} >> \Omega^1 (X;\ad\bar P)^-
@> d^+_{\bar A} + D_A \sigma >> \Omega^2_+ (X;\ad\bar P)^-.
\end{CD}
\]
\end{remark}

\begin{proposition}
There exist equivariant perturbations $\tau_k$ as in (\ref{E:tau-def'})
such that, for any small generic $\ep_k\in \mathbb R$ and $\tau = \sum\;
\ep_k \cdot \tau_k$, the moduli space $\M_{\tau}(P)$ is non--degenerate
at each of its eight--orbits.
\end{proposition}

\begin{proof}
First observe that the argument with Uhlenbeck's compactness theorem as in
the proof of Proposition \ref{P:gen} shows that, by choosing $\ep_k$ small
enough, one can make $\M^{\chi}_{\tau} (P)$ belong to an arbitrarily small
neighborhood of $\M^{\chi}(P)$.

Since $\M(P)$ is assumed to be non--degenerate in the tangential direction,
its local structure at any $[A] \in \M^{\chi}(P)$ is given by the Kuranishi
map
\[
D_A\tau = \sum\; \ep_k\cdot D_A\tau_k: H^1(X;\ad A)^- \to H^2_+(X;\ad A)^-.
\]
According to Lemma \ref{L:tau_k}, this map is an isomorphism for a generic
choice of $\ep_k$, hence its kernel vanishes at each point in the orbit of
$[A]$. If the $\ep_k$ are chosen to be sufficiently small, this guarantees
the non--degeneracy of $\M_{\tau}(P)$ at all nearby eight--orbits.

Repeating this argument for the rest of the (finitely many) $[A] \in
\M^{\chi}(P)$ and keeping in mind Remark \ref{R:pert}, one achieves
the non--degeneracy of $\M_{\tau} (P)$ at each of the eight--orbits.
\end{proof}


\section{The sixteen--orbits}\label{S:sixteen}
Throughout this section, by an equivariant admissible perturbation $\sigma$
we will always mean a linear combination of equivariant admissible
perturbations of the types (\ref{E:sigma'}) and (\ref{E:tau-def'}). Moreover,
we will assume that $\M_{\sigma}(P)$ is non--degenerate at the eight--orbits
and compact; that such $\sigma$ exists is the main result of the previous
section.

In this section, we will achieve the full non--degeneracy of $\M_{\sigma}(P)$
by further perturbing it at the 16--orbits. This task is simplified by the
fact that the stratum of 16--orbits, $\M^0_{\sigma}(P)\subset \M_{\sigma}(P)$,
is compact\,:  both four-- and eight--orbits in the compact $\M_{\sigma}(P)$
are non--degenerate and hence isolated. On the negative side, we will need to
perturb the already perturbed 16--orbits. These no longer consist of
projectively flat connections, hence the methods of Section \ref{S:tangent}
and specifically Lemma \ref{L:lifts} need to be adapted to this new situation.

\begin{lemma}\label{L:stab16}
Let $p: \tilde X \to X$ be the regular 16--fold covering of $X$ corresponding
to the mod 2 abelianization homomorphism $\pi_1 X \to H_1 (X; \Z_2)$.  If $A$
is a connection on $P$ whose stabilizer in $H^1 (X;\Z_2)$ is trivial then the
pull back connection $\tilde A$ is irreducible.
\end{lemma}

\begin{proof}
Let $x\in X$ and $\tilde x\in \tilde X$ be points such that $p(\tilde x) = x$,
then we have an exact sequence of monoids
\[
\begin{CD}
1 @>>> \Omega (\tilde X,\tilde x) @> p_* >> \Omega(X,x) @> q >> H_1(X;\Z_2)
@>>> 1
\end{CD}
\]
where $q(\gamma) = [\gamma]$ is the homology class of $\gamma$.  Write each
of the sixteen elements of $H_1(X;\Z_2)$ in the form $q(\gamma_k)$ for some
fixed loops $\gamma_k \in \Omega(X,x)$. Then every $h\in \im(\hol_A)\subset
SU(2)$ can be written in the form
\[
h = \hol_A(\gamma) = \hol_A (p_*(\tilde\gamma)\cdot\gamma_k) =
\pm\,\hol_{\tilde A}(\tilde\gamma)\cdot\hol_A(\gamma_k),
\]
where $\tilde\gamma\in\Omega(\tilde X,\tilde x)$. This means that
$\im(\hol_{\tilde A})$ has finite index in $\im(\hol_A)$.

Suppose that $\tilde A$ is reducible so that $\im (\hol_{\tilde A})$ is
contained in a copy of $U(1)$. Then $\im (\hol_A)$ is contained in a finite
extension of $U(1)$, that is, a copy of binary dihedral group $U(1)\cup j
\,U(1)$. Equivalently, if $\bar A$ is the connection on $\bar P$ adjoint to
$A$  then $\im (\hol_{\bar A}) \subset O(2)$ where the embedding $j: O(2)
\to SO(3)$ is given by $j(a) = \det(a) \oplus a$.  We wish to show that $A$
has a non--trivial stabilizer in $H^1 (X; \Z_2)$.

Let $\phi_{\a\b}:  U_\a \cap U_\b \to U(2)$ be a gluing cocycle for $P$ then
$\bar\phi_{\a\b}: U_\a \cap U_\b \to SO(3)$ is a gluing cocycle for $\bar P$.
Since $\bar P$ admits a connection $\bar A$ whose holonomy is $O(2)$ we may
assume that $\bar\phi_{\a\b}$ factors through the embedding $j: O(2) \to
SO(3)$. We wish to find a family of functions $u_\a: U_\a \to U(2)$
such that
\begin{equation}\label{E:patch}
\begin{aligned}
&u_\a \phi_{\a\b} = \phi_{\a\b} u_\b\ \ \text{on}\ \  U_\a \cap
U_\b,\ \text{and} \\
&u_\a du_\a^{-1} + u_\a A_\a u_\a^{-1} -
1/2\,h_\a\,dh_\a^{-1}\Id = A_\a \ \  \text{over}\ \
U_\a.
\end{aligned}
\end{equation}
Here, $h_\a = \det u_\a$.

The strategy is to make an initial choice $u_\a = i \in SU(2)$, and
then to modify that choice to make conditions~\eqref{E:patch} hold.
Since $\ad(i) = \diag(1,-1,-1)$ we have
$\ad(i)\bar\phi_{\a\b} = \bar\phi_{\a\b}\ad(i)$ in $SO(3)$. The same equality
holds in $U(2)$ modulo center, that is, $i\,\phi_{\a\b} = t_{\a\b}\phi_{\a\b}
\,i$ for some $t_{\a\b}: U_\a \cap U_\b \to U(1)$. Observe that $t_{\a\b}$ is
in the commutator subgroup of $U(2)$, hence, being central, only takes values
$\pm 1$.  Thus $t_{\a\b}$ defines an obstruction $t \in H^1 (X; \Z_2)$.

This obstruction need not vanish. On the other hand, we are allowed to redefine
$u_\a$ to be $t_\a i$ for some $t_\a: U_\a \to U(1)$. Because of this
extra freedom,
we only need vanishing of the image of $t$ in $H^1 (X; S^1) = H^2 (X; \Z)$.
This image equals $\beta(t) \in H^2 (X;\Z)$ where $\beta: H^1 (X;\Z_2) \to
H^2 (X;\Z)$ is the Bockstein operator. Since $H^2 (X; \Z)$ is torsion free,
$\beta(t) = 0$.

This means that there exist $t_\a: U_\a \to S^1$ such that  $u_\a = t_\a i$
satisfy the equation $u_\a\phi_{\a\b} = \phi_{\a\b}u_\b$ on $U_\a\cap U_\b$.
Finally, to verify the equation $u_\a du_\a^{-1} + u_\a A_\a u_\a^{-1} -
1/2\,h_\a\,dh_\a^{-1}\Id = A_\a$, we notice that $A$ is uniquely determined
by $\bar A$ and the determinant connection, hence it is enough to check the
above equality for connections $\bar A$ and $\tr A$. Both are immediate from
the definition of $u_\a$.
\end{proof}

\begin{lemma}\label{L:ab16}
For any $[A] \in \M^0_{\sigma}(P)$ there is a collection of mod-2 trivial
disjoint loops $\psi_k$ and a collection of self--dual $2$--forms  $\nu_k$,
each supported in its respective $\psi_k (S^1 \times D^3)$, such that the
forms $\Pi\hol_A (\psi_k)\, \otimes\,\nu_k$ generate the second cohomology
$H^2_{+,\sigma} (X; \ad A)$ of complex (\ref{D:pert}).
\end{lemma}

\begin{proof}
As in the proof of Lemma \ref{L:ab8}, we follow the proof of Lemma
2.5 in Donaldson \cite{donaldson:orientation}.

The subspace $H^2_{+,\sigma}(X;\ad A) = \coker(d^*_{\bar A} + D_A \sigma)$
of $\Omega^2_+ (X;\ad\bar P)$ is finite dimensional because of the
ellipticity of complex (\ref{D:pert}). Therefore, there exist finitely
many points $x_1,\ldots, x_m$ in $X$ such that the evaluation map
\[
\ev: H^2_{+,\sigma}(X;\ad A) \longrightarrow\; \bigoplus_{p = 1}^m\;
\Lambda^2_+ (\ad\bar P)_{x_p}
\]
given by $\ev(\omega) = (\omega(x_1),\ldots,\omega(x_m))$ is injective.
Note that, since $(\ad\bar P)_{x_p}$ is three--dimensional, each vector
space $\Lambda^2_+ (\ad\bar P)_{x_p}$ has dimension nine.

For each of the points $x_p$, consider the possible holonomies $\hol_A
(\psi)$ around mod-2 trivial loops $\psi$ based at $x_p$. The
restriction of $\hol_A$ onto mod-2 trivial loops coincides with the
map $\hol_{\tilde A}$ where $\tilde A$ is the pull back connection as
in Lemma \ref{L:stab16}.  According to that lemma, $\tilde A$ is
irreducible, hence one can find three mod-2 trivial loops $\psi_i$,
$i = 1, 2, 3$, such that the vectors $\Pi\hol_A (\psi_i)$ span
$(\ad\bar P)_{x_p}$.

Let $\nu_j$, $j = 1, 2, 3$, be self--dual $2$--forms supported near $x_p$
such that $\nu_j (x_p)$ span $(\Lambda^2_+)_{x_p}$. Then the nine vectors
$(\Pi\hol_A (\psi_i)\,\otimes\,\nu_j)(x_p)$ span the vector space
$\Lambda^2_+ (\ad\bar P)_{x_p}$. Repeating this construction for the rest
of the points $x_p$ results in a collection of loops and forms which we
will call $\psi_k$ and $\nu_k$, respectively. One may assume that $\supp
\nu_k \subset \psi_k(S^1\times D^3)$ for all $k$ and, after a small
isotopy, $\psi_k(S^1\times D^3)\,\cap\,\psi_l(S^1\times D^3) = \emptyset$
for all $k \ne l$. The result now follows because the supports of $\nu_k$
can be chosen so small that no non--zero form in $H^2_{+,\sigma}(X;\ad A)$
can be orthogonal to all of the $\Pi\hol_A(\psi_k)\,\otimes\,\nu_k$.
\end{proof}

The moduli space $\M^0_{\sigma}(P)$ is compact therefore we can suppose
that the mod-2 trivial loops $\psi_k$ and forms $\nu_k$ can be chosen
so that the conclusion of Lemma \ref{L:ab16} holds for all $[A] \in
\M^0_{\sigma}(P)$ simultaneously. We will fix such a choice.

\begin{proposition}\label{P:gen'}
There exists a real number $r > 0$ such that, with the choice of $\psi_k$
and $\nu_k$ as above and perturbation $\sigma'$ defined by the formula
(\ref{E:sigma'}), the moduli space $\M^0_{\sigma +\sigma'}(P)$ is
non--degenerate for a generic $\bar\ep = (\ep_1, \ldots, \ep_N)$ inside
the ball $\|\bar\ep\| < r$ in $\mathbb R^N$.
\end{proposition}

\begin{proof}
Given Lemma \ref{L:ab16}, the proof is almost identical to the proof of
Proposition \ref{P:gen}.
\end{proof}

Since small enough perturbations  $\sigma'$  used in Proposition \ref{P:gen'}
to make $\M^0_{\sigma} (P)$ non--degenerate do not spoil the non--degeneracy
at the four-- and eight--orbits, we are finished.


\section{Proof of Theorem \ref{T:main}}\label{S:main}
Let $P$ be an admissible $U(2)$--bundle on $X$, and $\bar P$ its associated
$SO(3)$--bundle.   If $w_2(\bar P)$ is not in the image of the map $\iota:
H^2 (\pi_1 X; \Z_2) \to H^2 (X;\Z_2)$, see (\ref{E:iota}), then the theorem
holds by Corollary \ref{C:iota}.

If $\M(P)$ is non--degenerate, consider the action on $\M(P)$ of $H^1(X;\Z_2)
= (\Z_2)^4$.  First, there are no orbits of orders one or two, see Remark
\ref{R:one-orbits}. According to Corollary \ref{C:det}, the contribution of
the four--orbits to $\bar\lambda (X,P)$ equals $\det X\pmod 2$. The theorem
now follows since the orbits consisting of eight and sixteen elements do not
contribute to $\bar\lambda(X,P)\pmod 2$ at all.

In general, $\M(P)$ becomes non--degenerate after an admissible equivariant
perturbation $\sigma$ as above. Since the four--orbits are
already non--degenerate,  they will remain such if $\sigma$ is sufficiently
small. The perturbation $\sigma$ will not create orbits with one or two
elements, or new orbits with four elements. Since $\sigma$ is equivariant,
the above argument discarding the orbits with more than four elements can be
applied again to show that $\bar\lambda (X,P) \equiv \det X\pmod 2$.

The statement of Theorem \ref{T:main} about the Rohlin invariant has been
verified in Corollary \ref{C:rhocup}.


\section{An example}\label{S:logtransform}

Theorem \ref{T:main} shows that, modulo 2, the invariant $\bar\lambda(X,P)$
equals the determinant of $X$.  It is natural to wonder if a similar result
holds over the integers. By analogy with the three--dimensional case, see
\cite{ruberman-saveliev:casson}, one would conjecture that $\bar\lambda(X,P)
=(\det X)^2$. In this section, we present a simple example that demonstrates
that this is not the case, and moreover, that $\bar\lambda(X,P)$ depends on
the choice of $P$.

The example may be described succinctly by viewing the $4$--torus as
an elliptic
fibration over $T^2$. Then we define a manifold $\Tqq$ to be the result of
log transforms of orders $q$ and $-q$ on two fibers.  For any value of $q$,
the manifold $\Tqq$ has the integral homology of $T^4$, and for $q$ odd
(which we assume to be the case henceforth) $\Tqq$ is an odd
$\Z[\Z]$--homology torus (in particular, it is spin). We will find that the
value of $\bar\lambda(\Tqq,P)$ depends on the choice of bundle $P$; for some
bundles it is given by $\pm q^2$ while for others it is equal to $\pm 1$.

This result is similar to the calculation of degree zero Donaldson
invariants of the log transform of an elliptic surface.  In fact the
manifold $\Tqq$ is the double branched cover (with exceptional curves
blown down) of a log transform of a K3 surface.  More generally, one
could consider homology tori coming from the double branched cover of
the Gompf-Mrowka non-complex homotopy K3 surfaces
\cite{gompf-mrowka}, which are log transforms of K3 surfaces along
two or three non-homologous tori.
The value of $\bar\lambda$ appears to be given by the analogue
of~\cite[Theorem 3.3]{gompf-mrowka} but the calculations are not
particularly revealing.  That the two results are similar is perhaps
not very surprising; compare~\cite{kronheimer:flat}.


\subsection{Log transforms on a torus}
Let us describe the construction in a little more detail.  Write $T^4
= S^1 \times S^1 \times S^1 \times S^1$ and let $u,v,x,y$ be the
coordinates on the four circles. Regard this torus as an elliptic
fibration via projection onto the first two coordinates, and call the
base torus $T$ and the fiber torus $F$.  Choose disjoint disks $D_a$,
$D_b$ in $T$, and let $\alpha$ and $\beta$ be their boundary curves,
oriented so that $\alpha + \beta = 0$ in $H_1(T - (D_a\cup D_b))$.
For all of these curves, we will use the same letter to denote the
curve and a coordinate on the curve.  Define
\[
\Tqq = ((T - (D_a\cup D_b)) \times F)\,\cup_{\phi_a}\,(D_a \times F)\,
\cup_{\phi_b}\,(D_b \times F)
\]
where $\phi_a: \p D_a \times F \to \p (T - (D_a\cup D_b)) \times F$ and
$\phi_b: \p D_b \times F \to \p (T - (D_a\cup D_b)) \times F$ are given
by the formulas
\begin{align*}
\phi_a(\alpha,x,y) & = (\alpha^q x^{q-1}, \alpha x, y)\quad \mathrm{and}\\
\phi_b(\beta,x,y) & = (\beta^{-q}x^{-q-1}, \beta x,y).
\end{align*}
It is not difficult to show that $\Tqq $ is a homology $4$--torus and that
its fundamental group $\pi_1(\Tqq)$ has presentation (with the above-named
curves representing fundamental group elements)
\begin{equation}\label{E:tqqgroup}
\langle x,y,u,v,\alpha,\beta \; |\; x,y\, \mathrm{central},\; \alpha^q
x= 1,\; \beta^{-q}x = 1,\; [u,v]\alpha\beta =1 \rangle.
\end{equation}
Note that projection onto the $y$-circle (the last $S^1$ factor) represents
$\Tqq$ as $\tqq \times S^1$, where $\tqq$ is obtained by Dehn filling of
$(T - (D_a\cup D_b)) \times S^1$.

The tori $u \times x$, $u \times y$, $v \times x$, $v \times y$,
$x \times y$ give generators for the second homology of $\Tqq$ with
$\Z_2$--coefficients (this is not true over the integers; as
in~\cite[\S2]{gompf-mrowka}, the regular fiber $F = x \times y$ is $q$
times the multiple fiber. But since we assume that $q$ is odd, we may as
well use the regular fiber). The remaining generator can be described as
an immersed surface in $\tqq$ as follows. Consider a disc $D$ in the base
torus $T$ containing $D_a$ and $D_b$. Start with $q$ copies of the twice
punctured torus $T - (D_a \cup D_b)$, whose boundary consists of $q$
parallel copies of $\alpha$ on one side,  and $q$ copies of $\beta$ on
the other. In the solid torus that is glued in to create the multiple
fiber (with meridian glued to $\alpha^q x$), there is a planar surface
with boundary the $q$ copies of $\alpha$ plus one copy of $x$; likewise
in the other solid torus, except with inverted orientation on $x$. The
two copies of $x$ are parallel in $D \times S^1$, so one can fill in an
annulus between them. The union of these is an immersed surface, say
$\Sigma\subset \tqq$, carrying the correct homology class.


\subsection{The pull-back case}\label{S:pull-back}
Let $P$ be a $U(2)$--bundle over $\Tqq$ that pulls back from a bundle over
$\tqq$ (which we call again $P$) having $w_2(\bar P) \ne 0$. Since $\Tqq$
is odd, the bundle $P$ is admissible. Let $\lppp(\tqq,w)$ with $w = w_2(P)$
be the Casson-type invariant of the homology $3$--torus $\tqq$ introduced in
\cite{ruberman-saveliev:casson}; we refer to that paper for all the
definitions. It is not difficult to show using techniques of
\cite{ruberman-saveliev:mappingtori} that the invariants $\bar\lambda
(\Tqq,P)$ and $\lppp(\tqq,w)$ coincide, at least when the
$SO(3)$--representation variety of $\pi_1(\tqq)$ is non--degenerate.

The remainder of this subsection will be devoted to a direct calculation of
the invariant $\lppp (\tqq,w)$ for two of the possible seven non--trivial
values of  $w \in H^2(\tqq;\Z_2)$ (a similar argument can be used to show
that the other five choices of $w$ give the same result, compare with
Remark \ref{R:comparison}). To describe our two choices for $w$, we present
the group $\pi_1(\tqq)$ as
\[
\langle x,u,v,\alpha,\beta\; |\; x\; \mathrm{central},\; \alpha^q x= 1,\;
\beta^{-q}x = 1,\; [u,v]\alpha\beta =1 \rangle,
\]
compare with (\ref{E:tqqgroup}), and consider representations $\rho:
\pi_1(\tqq)\to SO(3)$ whose $w_2 (\rho)$ is non--trivial on the generator
$x\times u \in H_2 (\tqq;\Z_2)$ and trivial on $x\times v$. This gives us
the two possibilities for $w = w_2 (\rho)$, depending on how it evaluates
on $\Sigma$, the third generator of $H_2(\tqq;\Z_2)$.

It is a standard fact \cite{braam-donaldson:knots} that since $w_2(\rho)$
evaluates non-trivially on the torus $x\times u$, the restriction of
$\rho$ to $x\times u$ is a $\ZZ$ representation, as described in Section
\ref{S:four-orbits}. Therefore, we may assume that $\rho(x)$ is a
$\pi$--rotation about the $x$--axis and $\rho(u)$ is a $\pi$--rotation
about the $y$--axis in $\mathbb R^3$. Note that since $\langle w_2(\rho),
x \times v\rangle = 0$, then $\rho(v)$ is a rotation around the $x$-axis
by an angle to be determined shortly.

It turns out that the representation $\rho$, and in particular the value
of $\rho(v)$, are completely determined by $\rho(\alpha)$, $\rho(\beta)$
and the evaluation of $w_2(\rho)$ on $\Sigma$. This is immediate from the
following recipe for computing $\langle w_2(\rho),\Sigma\rangle$. First,
choose lifts $\tilde\rho(x)$, $\tilde\rho(u)$, $\tilde\rho(v) \in SU(2)$
of the specified $SO(3)$ elements $\rho(x)$, $\rho(u)$ and $\rho(v)$. Of
course we know that $[\tilde\rho(x), \tilde\rho(u)] = -1$, because $w_2
(\rho)$ evaluates non-trivially on the torus $x\times u$. Now choose lifts
$\tilde\rho(\alpha)$ and $\tilde\rho(\beta)$ such that $[\tilde\rho(u),
\tilde\rho(v)]\tilde\rho(\alpha)\tilde\rho(\beta) = 1$. This can be done
by first choosing arbitrary lifts and then adjusting the value of (say)
$\tilde\rho(\beta)$. Then $\langle w_2(\rho),\Sigma \rangle$ is trivial
or non--trivial depending on whether $\tilde\rho(\alpha)^q
\tilde\rho(\beta)^q$ equals $1$ or $-1$.

To find $\rho(\alpha)$ and $\rho(\beta)$, we simply extract $q$-th roots of
$\rho(x)$ so that $\rho(\alpha)$ is a rotation about the $x$--axis through
$\pi k/q$ with odd $k = 1,\ldots ,2q-1$, and $\rho(\beta)$ is a rotation about
the $x$ axis through $\pi\ell/q$ with odd $\ell = 1,\ldots ,2q-1$. One can
conjugate the entire representation (using the $\pi$--rotation with respect
to the $y$--axis, for instance) without changing $\rho(x)$ and $\rho(u)$
but switching between $(k,\ell)$ and $(2q-k,2q-\ell)$. This leaves us with
pairs $(\rho(\alpha),\rho(\beta))$ no two of which are conjugate to each
other as long as $1 \le k \le q$. Among these, there is a special pair with
$k = \ell = q$, and a total of $(q^2 - 1)/2$ other pairs.

The special pair gives a unique $SO(3)$ representation with image in $\ZZ$
and given $w_2$.   This can be easily seen from the last relation, which in
this case says that $\rho(u)$ and $\rho(v)$ commute -- therefore, $\rho(v)$
is a rotation through $0$ or $\pi$ about the $x$-axis.  From the recipe, it
is immediate that this rotation angle is $0$ if and only if
$\langle w_2(\rho),x \times v\rangle = 0$. Once $w_2$ is fixed, this
$SO(3)$ representation lifts to the unique four--orbit of projective
$SU(2)$ representations; in particular, it is non--degenerate, see
Proposition 4.4 of \cite{ruberman-saveliev:casson}.

For each of the remaining pairs $(\rho(\alpha),\rho(\beta))$, we use the
last relation to find the rotation angle of $\rho(v)$. A little calculation
shows that this equation determines two possible angles, $\theta$ and
$\theta + \pi$.  As before, the evaluation $\langle w_2(\rho),x \times v
\rangle$ is $0$ for exactly one of these choices. Therefore, for a given
choice of $w_2$, there is  a unique representation $\rho$ for each of the
pairs $(\rho(\alpha),\rho(\beta))$. It lifts to an eight--orbit of
projective $SU(2)$ representations -- this follows from the observation
that the image of $\rho$ belongs to $O(2) \subset SO(3)$ but not to a copy
of $\ZZ$. All these representations are non--degenerate (proving this is
an easy exercise in linear algebra) hence each of them is counted in
$\lppp(\tqq,w)$ with multiplicity two.

Thus we have a total of $2 \cdot (q^2-1)/2 + 1 = q^2$ representations
(counted with multiplicities). All of them should be counted with the same
sign -- for instance, because they all give rise to holomorphic bundles so
the orientations are the same at all points. Therefore, $\lppp(\tqq,w) =
\pm\,q^2$.

Another way to prove the latter formula (without referring to holomorphic
bundles) is by using surgery theory. A framed link describing $\tqq$ can
be obtained from the Borromean rings, with all three components framed by
zero, by adding two meridians of one of the components with respective
framings $q$ and $-q$. Surgering out one component of the Borromean rings
at a time and using Casson's surgery formula, one reduces the calculation
to that of Casson's invariant of a certain plumbed manifold. The latter is
given by an explicit formula. We leave details to the reader. The downside
of this surgery approach is that it does not give any information about
the non--degeneracy of the representation variety, which we use to identify
$\bar\lambda(\Tqq,P)$ with $\lppp(\tqq,w)$.

\begin{remark}\label{R:comparison}
Conjecturally, the gauge-theoretic invariant $\lppp(Y,w)$ of a homology
$3$--torus $Y$ is independent of $w$ and coincides (up to an overall sign)
with the combinatorially defined invariant of Lescop \cite{lescop:casson}.
According to \cite{ruberman-saveliev:casson} this conjecture holds modulo
2, and one hopes that it can also be verified over the integers by either
establishing a surgery formula in the gauge theory context, or by an
argument analogous to Taubes' theorem~\cite{taubes:casson,masataka,clm:III}.
Lescop shows that her invariant is given by $(\det Y)^2$, which is manifestly
independent of $w$. In the case at hand, $Y = \tqq$, the determinant equals
$q$ so the above calculation of $\lppp(\tqq,w) = \pm q^2$ weighs in for the
conjecture.
\end{remark}


\subsection{The non-pull-back case}
The case when the bundle $P$ does not pull back from $\tqq$ yields a different
answer for $\bar\lambda(\Tqq,P)$. Suppose that $\rho$ is an $\SO(3)$
representation with $\langle w_2(\rho), x\times y\rangle \neq 0$, and assume
that $\rho(x)$ is a $\pi$--rotation about the $x$-axis and $\rho(y)$ is a
$\pi$--rotation about the $y$-axis. Note that, since any element that
centralizes a $\ZZ$ subgroup in $SO(3)$ is an element of that subgroup,
we conclude that $\rho(u)$, $\rho(v)$, $\rho(\alpha)$ and $\rho(\beta)$
all belong to the same $\ZZ$ subgroup as $\rho(x)$ and $\rho(y)$. Since
$\rho(\alpha)$ is a $q$-th root of the $\pi$-rotation about the $x$-axis,
it also is a $\pi$-rotation about the $x$-axis. Thus $\rho(\alpha) =
\rho(x)$, and similarly, $\rho(\beta) = \rho(x)$.

Now, the relations in (\ref{E:tqqgroup}) hold for arbitrary choices of
$\rho(u), \rho(v) \in \ZZ$ --  in fact, fixing $\rho(u)$ and $\rho(v)$
amounts to fixing a particular value for $w_2(\rho)$. For instance, if
$w_2(\rho)$ vanishes on the $2$--torus $x \times u$, then $\rho(u)$ must
be either identity or a $\pi$--rotation about the $x$--axis. On the
other hand, if $\langle w_2(\rho),x \times u \rangle \neq 0$, then
$\rho(u)$ is a $\pi$--rotation around the $y$--axis or the $z$--axis.
Thus our freedom in choosing $\rho(u)$ and $\rho(v)$ is only restricted
by the requirement that $w_2 (\rho)\cup w_2 (\rho)$ be zero modulo 4.

Choose for instance $w_2 (\rho) = x^*\cup y^*$ where $x,y\in H^1(\Tqq;
\Z_2)$ are dual to the curves $x$ and $y$.  Then obviously $w_2 (\rho)
\cup w_2 (\rho) = x^*\cup y^*\cup x^*\cup y^* = 0\pmod 4$, and one can
easily see that this choice corresponds to having $\rho(u) = \rho(v) =
1$. For this particular choice of $w_2$, the above $\ZZ$ representation
$\rho$ is the only point in the representation variety. By Proposition
\ref{P:4-nondeg}, such a representation is automatically non--degenerate
(remember that since $q$ is odd, $\Tqq$ is an odd homology torus), so the
invariant equals $\pm 1$.


\newpage

\begin{thebibliography}


\bibitem{braam-donaldson:knots} \textbf{P\,J Braam}, \textbf{S\,K
Donaldson}, \emph{Floer's work on instanton homology, knots and
surgery}, from: ``The Floer memorial volume'', Progr.  Math. 133,
Birkh\"auser, Basel (1995) 195--256 \MR{1362830}

\bibitem{brown:cohomology}
\textbf{K\,S Brown}, \emph{Cohomology of groups}, Graduate Texts in Mathematics
  87, Springer-Verlag, New York (1994) \MR{1324339}

\bibitem{clm:III}
\textbf{S\,E Cappell}, \textbf{R Lee}, \textbf{E\,Y Miller}, \emph{Self-adjoint
  elliptic operators and manifold decompositions. {III}. {D}eterminant line
  bundles and {L}agrangian intersection}, Comm. Pure Appl. Math. 52 (1999)
  543--611 \MR{1670052}

\bibitem{donaldson}
\textbf{S\,K Donaldson}, \emph{An application of gauge theory to
  four-dimensional topology}, J. Differential Geom. 18 (1983) 279--315
  \MR{710056}

\bibitem{donaldson:orientation}
\textbf{S\,K Donaldson}, \emph{The orientation of {Y}ang-{M}ills moduli spaces
  and {$4$}-manifold topology}, J. Differential Geom. 26 (1987) 397--428
  \MR{910015}

\bibitem{donaldson:floer}
\textbf{S\,K Donaldson}, \emph{Floer homology groups in {Y}ang-{M}ills theory},
  Cambridge Tracts in Mathematics 147, Cambridge University Press, Cambridge
  (2002) \MR{1883043}

\bibitem{floer:instanton}
\textbf{A Floer}, \emph{An instanton-invariant for {$3$}-manifolds}, Comm.
  Math. Phys. 118 (1988) 215--240 \MR{956166}

\bibitem{furuta:perturbation}
\textbf{M Furuta}, \emph{Perturbation of moduli spaces of self-dual
connections}, J. Fac. Sci. Univ. Tokyo Sect. IA Math. 34 (1987) 
275--297

\bibitem{furuta-ohta}
\textbf{M Furuta}, \textbf{H Ohta}, \emph{Differentiable structures on
  punctured {$4$}-manifolds}, Topology Appl. 51 (1993) 291--301 \MR{1237394}

\bibitem{gompf-mrowka}
\textbf{R\,E Gompf}, \textbf{T\,S Mrowka}, \emph{Irreducible {$4$}-manifolds
  need not be complex}, Ann. of Math. (2) 138 (1993) 61--111 \MR{1230927}

\bibitem{herald:perturbations}
\textbf{C Herald}, 
\emph{{Legendrian cobordism and Chern--Simons theory on
$3$--manifolds with boundary}}, Comm.\ Anal.\ Geom. 2 (1994) 
337--413 \MR{1305710}

\bibitem{milnor-husemoller}
\textbf{D Husemoller}, 
\textbf{J Milnor}, \emph{Symmetric Bilinear Forms}, 
Ergebnisse series 73, Springer-Verlag, New York--Heidelberg (1973)
\MR{0506372}

\bibitem{kaplan:even}
\textbf{S\,J Kaplan}, \emph{Constructing framed {$4$}-manifolds with given
  almost framed boundaries}, Trans. Amer. Math. Soc. 254 (1979) 237--263
  \MR{539917}

\bibitem{kronheimer:higher-rank}
\textbf{P\,B Kronheimer}, 
\emph{Four--manifolds invariants from higher-rank
  bundles} (2004), 
\url{http://www.math.harvard.edu/~kronheim/higherrank.pdf}

\bibitem{kronheimer:flat}
\textbf{P\,B Kronheimer}, \emph{Instanton invariants and flat connections on the
{Kummer} surface}, Duke Math. J. 64 (1991) 
229--241 \MR{1136374}

\bibitem{laitinen:ends}
\textbf{E Laitinen}, \emph{End homology and duality}, Forum Math. 8 (1996)
  121--133 \MR{1366538}

\bibitem{lescop:casson} \textbf{Christine Lescop}, \emph{Global
surgery formula for the {C}asson--{W}alker invariant}, Annals of
Mathematics Studies 140, Princeton University Press, Princeton, NJ
(1996) \MR{1372947}

\bibitem{masataka}
\textbf{K Masataka}, \emph{Casson's knot invariant and gauge theory}, Topology
  Appl. 112 (2001) 111--135 \MR{1823600}

\bibitem{milnor:covering}
\textbf{J\,W Milnor}, \emph{Infinite cyclic coverings}, from: ``Conference on
  the Topology of Manifolds (Michigan State Univ., E. Lansing, Mich., 1967)'',
  Prindle, Weber \& Schmidt, Boston, Mass. (1968)  115--133 \MR{0242163}

\bibitem{pontrjagin:pi3} \textbf{L Pontrjagin}, \emph{Mappings of the
three-dimensional sphere into an {$n$}-dimensional complex},
C. R. (Doklady) Acad. Sci. URSS (N. S.)  34 (1942) 35--37 \MR{0008135}

\bibitem{ruberman:ds}
\textbf{D Ruberman}, \emph{Doubly slice knots and the {C}asson-{G}ordon
  invariants}, Trans. Amer. Math. Soc. 279 (1983) 569--588 \MR{709569}

\bibitem{ruberman-saveliev:casson}
\textbf{D Ruberman}, \textbf{N Saveliev}, \emph{Rohlin's invariant and gauge
  theory. {I}. {H}omology 3-tori}, Comment. Math. Helv. 79 (2004) 618--646
  \MR{MR2081729}

\bibitem{ruberman-saveliev:mappingtori}
\textbf{D Ruberman}, \textbf{N Saveliev}, \emph{Rohlin's invariant and gauge
  theory. {II}. {M}apping tori}, \gtref{8}{2004}{2}{35}{76} \MR{MR2033479}

\bibitem{ruberman-saveliev:survey}
\textbf{D Ruberman}, \textbf{N Saveliev}, \emph{{C}asson--type
invariants in dimension four}, from: ``{G}eometry and {T}opology of
{M}anifolds'', Fields Institute Communications 47, AMS (2005), 
\arxiv{math.GT/0501090}

\bibitem{ruberman-strle:tori}
\textbf{D Ruberman}, \textbf{S Strle}, \emph{Mod 2
{S}eiberg-{W}itten invariants of homology tori}, Math. Res. Lett. 7 
(2000) 789--799
\MR{1809302}

\bibitem{taubes:periodic}
\textbf{C\,H Taubes}, \emph{Gauge theory on asymptotically periodic
  {$4$}-manifolds}, J. Differential Geom. 25 (1987) 363--430 \MR{882829}

\bibitem{taubes:casson} \textbf{C\,H Taubes}, \emph{Casson's invariant
and gauge theory}, J. Differential Geom.  31 (1990) 547--599
\MR{1037415}

\bibitem{turaev:linking}
\textbf{V\,G Turaev}, \emph{Cohomology rings, linking coefficient forms and
  invariants of spin structures in three-dimensional manifolds}, Mat. Sb.
  (N.S.) 120(162) (1983) 68--83, 143 \MR{687337}

\bibitem{uhlenbeck:lp-bounds}
\textbf{K Uhlenbeck}, \emph{Connections with {$L^p$}-bounds on curvature},
Comm.\ Math.\ Phys. 83 (1982) \MR{0648356}
31--42

\bibitem{witten:monopole}
\textbf{E Witten}, \emph{Monopoles and four-manifolds}, Math. Res. Lett. 1
  (1994) 769--796 \MR{1306021}

\end{thebibliography}
\end{document}